\documentclass[twoside]{article}
\usepackage{arxiv}
\usepackage{fancyhdr}
\usepackage[font=small]{caption}
\usepackage{float}
\usepackage[utf8]{inputenc} 
\usepackage[T1]{fontenc}    
\usepackage{hyperref}       
\usepackage{url}            
\usepackage{pdfpages}
\usepackage{booktabs}       
\usepackage{biblatex} 
\usepackage{amsfonts,amsmath}       
\usepackage{nicefrac}       
\usepackage{microtype}      
\usepackage{lipsum}		
\usepackage{graphicx}
\usepackage{doi}
\usepackage{caption}
\usepackage{subcaption}
\usepackage{ulem} 
\usepackage{caption}
\usepackage{subcaption}
\usepackage{microtype}
\usepackage{todonotes}
\usepackage[noend]{algpseudocode}
\usepackage{algorithm}
\usepackage{placeins}

\usepackage{fancyhdr}
\title{Data-driven macroscopic dynamics of complex networks using Topological Data Analysis and the Equation-Free Method}

\author{ \href{https://orcid.org/0000-0002-5651-5631}{\includegraphics[scale=0.06]{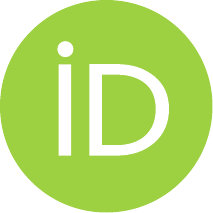}\hspace{1mm}Konstantinos Spiliotis}\thanks{Corresponding author.} \\
	Department of Civil Engineering, School of Engineering\\ Democritus University of Thrace, Xanthi, Greece\\
	\texttt{kspiliot@civil.duth.gr} \\
	\And
	\href{https://orcid.org/0000-0000-0000-0000}{\includegraphics[scale=0.06]{orcid.pdf}\hspace{1mm}Ole Sönnerborn} \\
	Department of Mathematics and Computer Science\\
    Karlstad University, Universitetsgatan 2\\ SE-651 88 Karlstad, Sweden
	\texttt{Ole.Sonnerborn@kau.se} 
    \And
	\href{https://orcid.org/0000-0002-1726-4892}{\includegraphics[scale=0.06]{orcid.pdf}\hspace{1mm}Haralampos Hatzikirou} \\
	Mathematics Department\\Khalifa University PO Box 127788\\Abu Dhabi, United Arab Emirates
	\texttt{haralampos.hatzikirou@ku.ac.ae} 
    \And
	\href{https://orcid.org/0000-0002-9743-8636}{\includegraphics[scale=0.06]{orcid.pdf}\hspace{1mm}Nikos I. Kavallaris} \\
	Department of Mathematics and Computer Science\\
    Karlstad University, Universitetsgatan 2\\ SE-651 88 Karlstad, Sweden
	\texttt{Nikos.Kavallaris@kau.se} 
}

\renewcommand{\shorttitle}{\textit{arXiv} Template}

\hypersetup{
pdftitle={Data-driven macroscopic dynamics of complex networks using Topological Data Analysis and the Equation-Free Method},
pdfsubject={q-bio.NC, q-bio.QM},
pdfauthor={David S.~Hippocampus, Elias D.~Striatum},
pdfkeywords={First keyword, Second keyword, More},
}

\begin{document}
\maketitle

\begin{abstract}
	In this work, we present a computational framework for exploring and analyzing the macroscopic dynamics of complex agent-based network models by integrating Topological Data Analysis with the Equation-Free Method. To demonstrate the effectiveness of our method, we apply it to Erdős--Rényi-type random networks.
Central to our approach is a Topological Data Analysis-based filtration process driven by the density of activated network nodes (agents), from which we extract a coarse-grained macroscopic topological observable. This observable is defined via persistent Betti numbers, thus requiring significantly reduced data dimensionality while retaining essential topological features. Subsequently, within the Equation-Free Method framework, we show firstly that a \textit{lifting procedure} can be achieved using topological properties and secondly, a data-driven evolution law that governs the dynamics of this macroscopic variable. Finally, we perform a numerical bifurcation and stability analysis to investigate the global behavior and qualitative transitions of the emergent macroscopic dynamics.
\end{abstract}


\section{Introduction}
{C}oupled dynamical systems arise in a wide range of physical and biological phenomena. These systems often exhibit multiple temporal and spatial scales, resulting from the interactive behavior of their constituent components. Such components may include particles in physical systems, cells in biological systems, or individual humans in the context of a pandemic---each contributing to the emergence of complex networks behavior. \cite{Str02, Kerch05, Bick20}.
%
%

Characteristic examples of complex networks can be found in systems neuroscience
, where millions of interacting neurons give rise to emergent behavior. The dynamics of individual neurons, modeled as nonlinear systems, were successfully captured by the seminal Hodgkin--Huxley model \cite{Hod52}, which earned the authors the 1963 Nobel Prize in Physiology (although the model was mathematical). This formalism continues to serve as a foundational framework for modeling large-scale brain networks \cite{Hod52, Ter02, Spil21}. Within this context, the structure of network connectivity and the architecture of coupling are critical in determining the dynamical properties of the system. Prominent phenomena emerging from these networks include self-organized criticality \cite{Spil11, Deco08, Deco13}, sustained oscillatory dynamics, and the propagation of traveling waves \cite{Lai02, Palk23, Bha22}.

Topological Data Analysis (TDA), a branch of data science, focuses on identifying patterns and simplifying the ``shape'' of data.
Rooted in algebraic topology, TDA reveals intrinsic topological features of data, such as voids and cavities, thereby offering powerful tools for uncovering complex structures and patterns \cite{Sag18, Chaud19,Volkening20}. These insights contribute to the development of machine learning methodologies by enabling the identification of structural features that are essential to the data's organization. In particular, techniques such as persistent homology \cite{Ada09} provide a means to encode both local geometric and global topological information into vectorized representations, which can be seamlessly integrated into machine learning pipelines \cite{Adams21}.


TDA has found increasing application in the study of complex dynamical systems, including agent-based models that simulate collective behaviors such as bird flocks, fish schools, and insect swarms \cite{Xian22,Ghar24,Corc17}. Beyond structural characterization, TDA contributes to the analysis of system dynamics and the prediction of critical phase transitions \cite{Male16, xian2021capturing, topa15, Cole21,Volkening20}. For instance, the study in Ref.~\cite{topa15} demonstrates that topological invariants such as Betti numbers can effectively quantify macroscopic behaviors of aggregates, complementing traditional metrics like polarization and angular momentum. Similarly, in Ref.~\cite{Cole21}, Betti numbers are employed to identify phase transitions in large-scale systems, drawing analogies to macroscopic properties such as magnetization in Ising-like models.

The models discussed above, whether agent-based aggregate models or Ising-type systems, share a fundamental characteristic: Although their macroscopic behavior is of primary interest, the governing evolution laws are defined exclusively at the microscopic level. The  Equation-Free Method (EFM), introduced in a series of studies \cite{Kev09, Gear03, Gear05, SIET11, Armaou04, Marsch14, Sieb17}, provides a computational framework for multiscale analysis of complex systems in the absence of explicit macroscopic equations. EFM enables tasks such as bifurcation analysis and stability assessment by leveraging short bursts of appropriately initialized microscopic simulations, without requiring closed-form macroscopic models \cite{Kev09}. This approach justifies its designation as ``equation-free'', as it facilitates system-level analysis even when macroscopic evolution equations are conceptually valid but not explicitly available. Suitable for complex systems where macroscopic equations are far from trivial to derive, the EFM has become a well-established method for analyzing such systems \cite{Rajen11, Bold14, Goun16, Bartal17, Rajendran17, MATTHEWS2019, Choi20161165}. It has been applied across a wide range of domains, from neuronal networks \cite{LAING08, Spil11, Spil24f, Kemeth18} to epidemic and social networks \cite{TSOU10, SIET11}, as well as crowd and traffic dynamics \cite{Marschel14DDif, Sieb18}. More recently, EFM has also been incorporated into machine learning approaches for complex systems \cite{Evan24, gallos21, PATSATZIS2023}.

In this study, we integrate TDA with the EFM to investigate and characterize the macroscopic dynamics of a complex random network of Erdős--Rényi type \cite{Koz05, Spil11, Spil10}. As a first step, the dynamic activity of individual agents, such as neurons, is projected onto the unit circle. Following this, we employ a filtration process using witness simplicial complexes to compute the minimal filtration radius at which the unit circle becomes topologically filled, corresponding to the emergence of the first Betti number, $\mathrm{Betti}_1$ \cite{Carl09}. This minimal filtration radius serves as a topological descriptor that correlates with the density of activated nodes in the network.

\begin{figure*}[t!]
\begin{center}
\hspace{-1.5cm}%
\hspace*{0.5cm} \includegraphics[width=0.7 \linewidth]{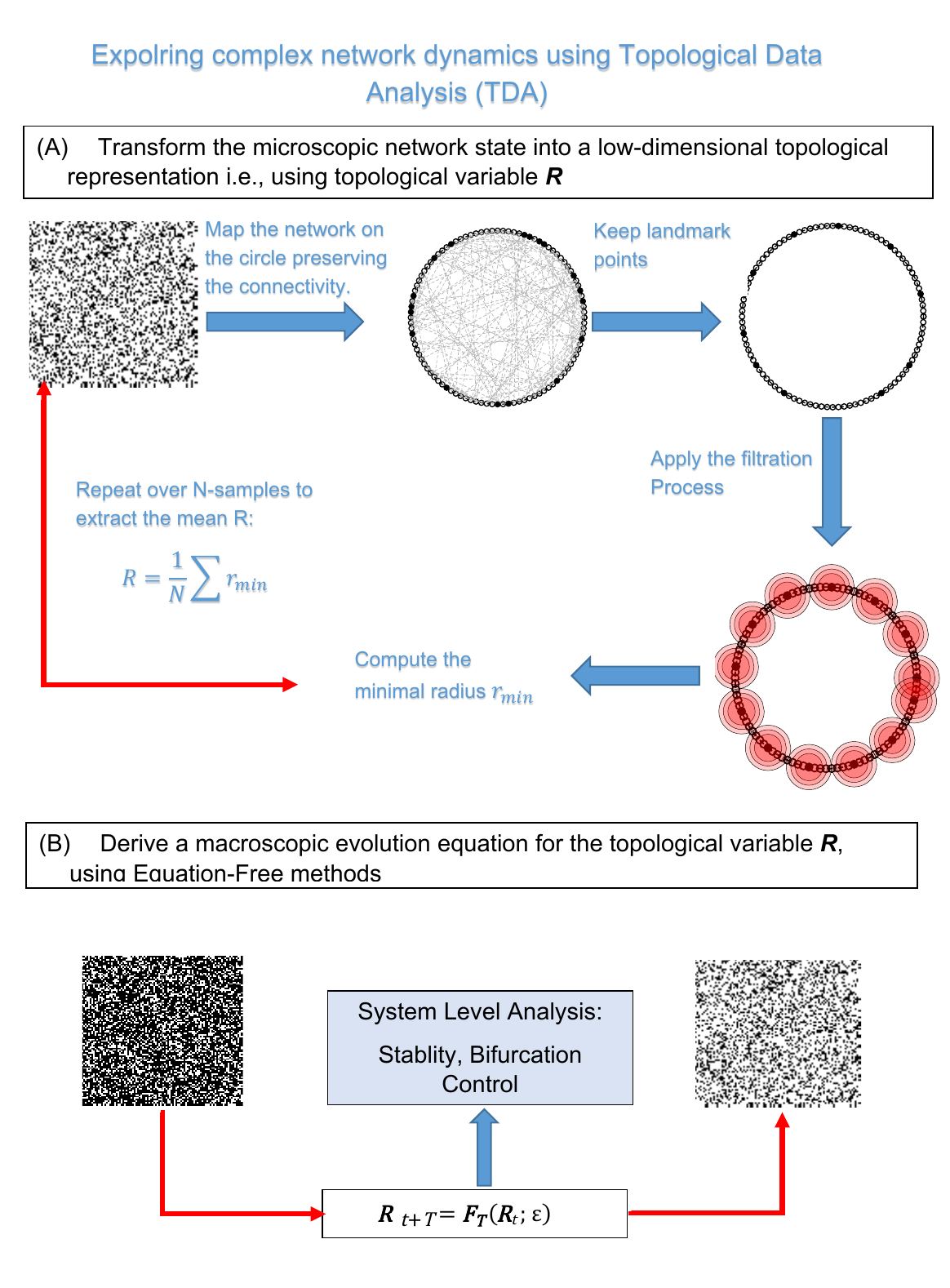}
\end{center}
\caption{Pipeline of the proposed method to analyse the network dynamics, combining topological data analysis with Equation-Free method \cite{Gear03,Kev09,Spil11,Marsch14}. \textbf{(A)} The microscopic network state can be transformed into a low-dimensional topological representation. To achieve this, we project the network state into the circle (where each agent in the network is represented now as a node in the unit circle), keeping the connectivity structure. Then, we select a small number of landmark points from the set of nodes in the circle. Using the filtration process \cite{Adams21,topa15}, we compute the minimum radius such that Betti1=1. Since our system is stochastic (at each time step), we repeat the same procedure and we average to extract the macroscopic topological value $R$. This procedure is described in section \ref{sec:topology} \textbf{(B)} Since we have a low dimensional representation of the network dynamics, i.e., the topological variable $R$, our next objective is to establish (numerically) an evolutionary equation: $R_{t+1}= F_T(R_t;\epsilon)$, using the fine atomistic rules of the agent-based model. This is implemented under the Equation-Free method \cite{Gear03,Kev09} and presented in section \ref{sec:results}. At this point, one can wrap around the coarse time stepper $F_T$ numerical solvers, e.g., Newton or Jacobian-free methods, to compute equilibria or periodic solutions (even unstable ones).} 
\label{fig:pip}  
\end{figure*}

During the filtration process, it is expected that densely activated neurons, when projected onto the unit circle, require a smaller radius to complete the topological filling of the circle compared to configurations with lower activation density. This observation enables a direct link between macroscopic dynamics and the topological feature known as the $\mathrm{Betti}_1$ radius, using a significantly reduced dataset comprising only the activated neurons. This connection is established through the implementation of a lifting operator, introduced in our novel approach (see Section~3.\ref{sec:lifting_operator_1}), which enhances previous methods based on simulated annealing by leveraging geometric information obtained from TDA to reconstruct a consistent microscopic network state from the topological signature of the minimal filtration radius. Within the EFM framework, we then analyze the macroscopic network dynamics as a function of this topological invariant. Additionally, we perform numerical bifurcation and stability analyses of the resulting macroscopic behavior. An overview of the proposed methodology is presented in Fig.~\ref{fig:pip}.

\section{An agent-based model with discrete dynamics}
\label{sec:model}

We consider a discrete state agent-based model. 
Agent-based models contain many interacting agents (e.g., cells, neurons, individuals), 
and the emergent macroscopic dynamics is far from trivial to predict \cite{PATSA23, SIET11, WANG15}. 
In a discrete state agent-based model, each state of an agent is typically described by an integer between $0$ and $n$, 
with the values between $1$ and $n-1$ corresponding to states within the refractory 
or recovery period. If the state is $n$, the agent is activated, and $0$ corresponds to the rest state. In our model, we only use two states corresponding to an inactive (rest) and an active (firing) state, {as in, for example, a neuronal network model with a refractory period or a susceptible-infected epidemiological network.}


\subsection{Description of the network}
\label{sec: Description of the network}
We consider a random network $(G, E)$, where $G$ is the set of agents (nodes) and $E$ is the set of edges, representing interactions between the agents.
We enumerate the agents using integers $1$ to $K$, with $K$ being the total number of agents.
In the network, any two agents $i$ and $j$ are connected by an edge $[i, j]$ with connectivity probability $p$.
We define the set of neighbors of agent $i$ to be the set of agents connected to it by an edge, $G_i=\{ j\in G : [i,j]\in E\}$. 
Furthermore, as mentioned above, we assume each agent has two states, $0$ (inactive) and $1$ (active), and we write $s_{i;t}$ for the state of the $i$th agent at time $t$. 
The density of active agents at time $t$ is then defined as
\begin{equation}
    d_t=\frac{1}{K}\sum_{i\in G} s_{i;t}.
\end{equation}
This density serves as a macroscopic variable (order parameter) that characterizes the network's state; 
a value of $d_t$ close to $1$ indicates that a significant portion of the agents are active, 
whereas a value close to $0$ indicates that the vast majority of agents are inactive.


The dynamics of the activity of the network is defined as follows. We start with a given activity density $d_0$ and randomly activate $K\cdot d_0$ of the agents. ($K\cdot d_0$ should be rounded off to an integer.) Let the activation probability $\epsilon$ be any number between, but not including, $0$ and $0.5$.
At each computational step we update the state of each agent according to the following majority rule:
\begin{itemize}
    \item Suppose the $i$th agent is inactive.
    Then, the $i$th agent is activated with probability $\epsilon$ if
    at most half of its neighbors are active, and it is activated with probability $1-\epsilon$
    if more than half of its neighbors are active.
    \item Suppose the $i$th agent is active.
    Then, the $i$th agent is deactivated with probability $\epsilon$ if
    at most half of its neighbors are inactive, and it is deactivated with probability $1-\epsilon$
    if more than half of its neighbors are inactive.
\end{itemize}

\subsection{Network dynamics}
In this work, we explore how large-scale collective behavior emerges from the network dynamics governed by the majority rule. Our analysis focuses on nonlinear phenomena, including the existence of multiple macroscopic stable phases, and phase transitions that occur within the network \cite{BOCCA16, Spil11, Spil24f}.

For the Erd\"os-R\'enyi type network described above, we have a fairly detailed picture of phase transitions \cite{Spil10,Spil11}. In particular, 
for a connectivity probability $p$ above a critical value $p_c$, the network has only one stable phase, independent of the value of the activation probability $\epsilon$. However, for $p<p_c$, there is a critical value $\epsilon_c$ for the activation probability below which the network has two locally stable collective states and goes through a phase transition, whereas only one stable state occurs for $\epsilon>\epsilon_c$; see \cite{Koz05, Spil11}. Figure \ref{fig:temp_kozma} shows the time evolution of the activity density for connectivity probability $p = 0.001$. The system shows two stable phases for small values of $\epsilon$: one corresponds to a low-activity phase, the other to a high-activity phase (Fig.~\ref{fig:temp_kozma}, plots (\textbf{a})--(\textbf{c})). The high-density phase disappears if we let $\epsilon$ pass $\epsilon_c$, and only one stable low-density phase exists (Fig.~\ref{fig:temp_kozma}, plot (\textbf{d})). This low-density phase remains for higher values of $\epsilon$. See the Supplementary Material for a detailed analysis. {In this study, we propose an alternative method based on persistent homology to analyze the system’s dynamics—including stability, bifurcation behavior, and phase transitions—by exploiting the topological properties of the underlying dynamics and the presence of \textit{landmark} points, which represent a significantly smaller subset of agents used to capture the system’s behavior.

The stability and bifurcation behavior of Erd\"os--R\'enyi type networks has been extensively investigated using the active agent density $d_t$ as a key metric \cite{Spil11}. In this study, we present a framework that integrates TDA and EFM, offering a highly computationally efficient method for analyzing stability and bifurcation phenomena in such networks. Furthermore, the proposed approach is generalizable and can be readily applied to a broad class of agent-based network models.

\begin{figure}[t]
    \centering
    \begin{subfigure}{.44\textwidth}{\includegraphics[width=1.1\textwidth]{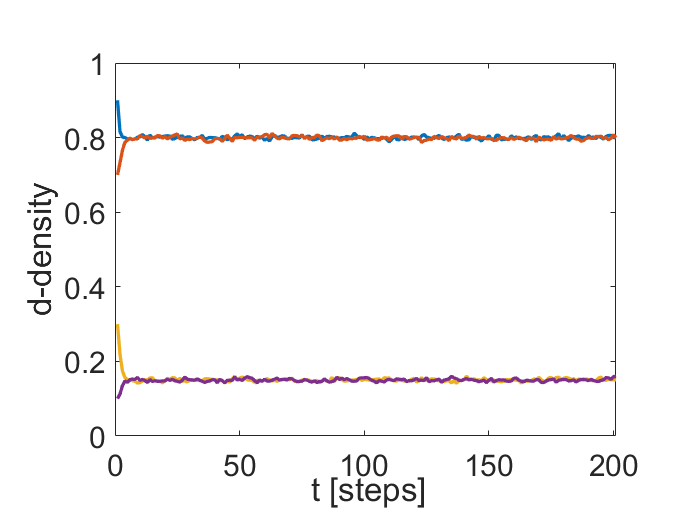}} 
    \caption{}
    \label{fig:No_rain_deg}
    \end{subfigure}
    \begin{subfigure}{.44\textwidth}{\includegraphics[width=1.1\textwidth]{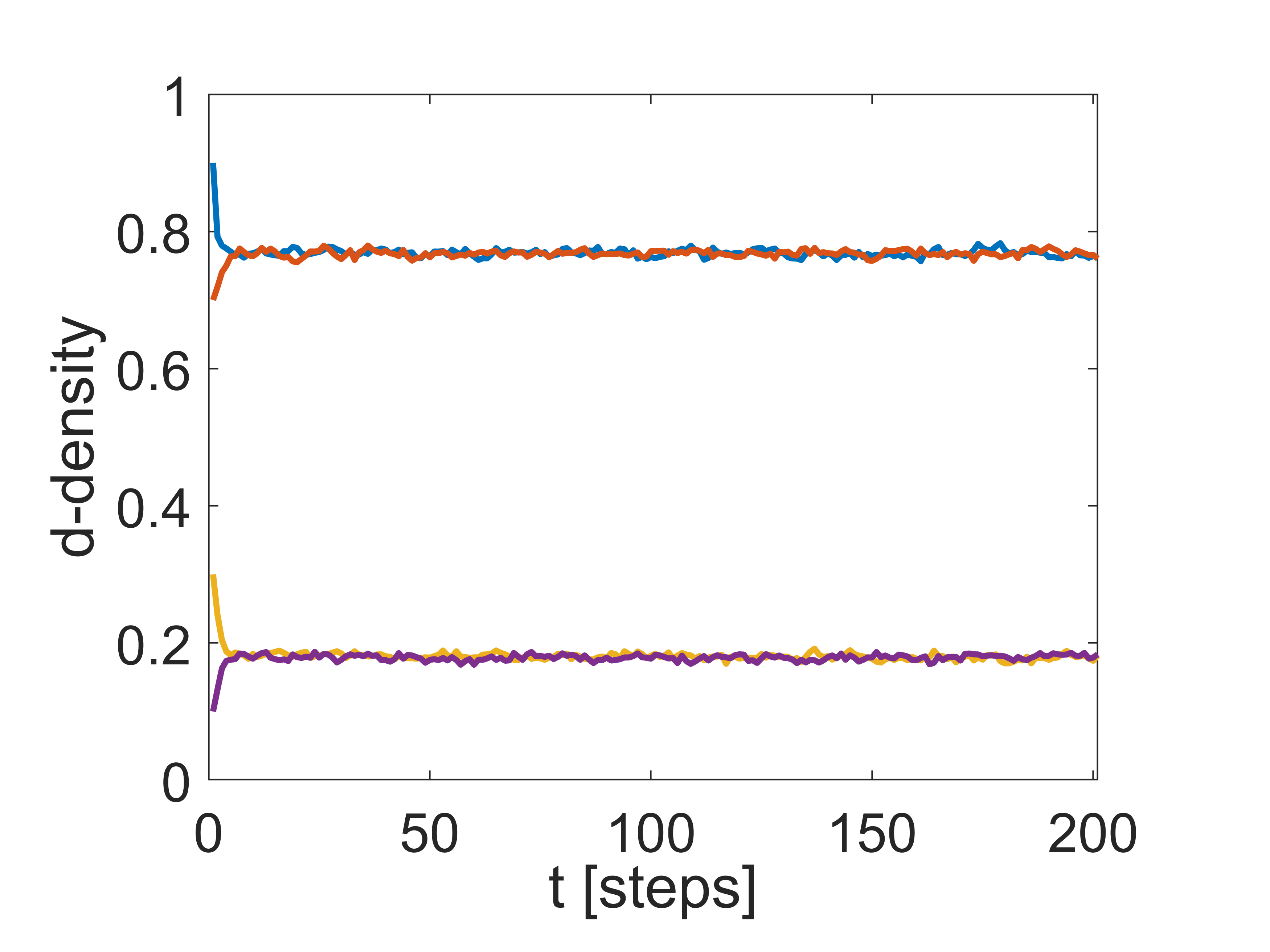}} 
    \caption{}
    \label{fig:No_rain_BC}
    \end{subfigure}
    \begin{subfigure}{.44\textwidth}{\includegraphics[width=1.1\textwidth]{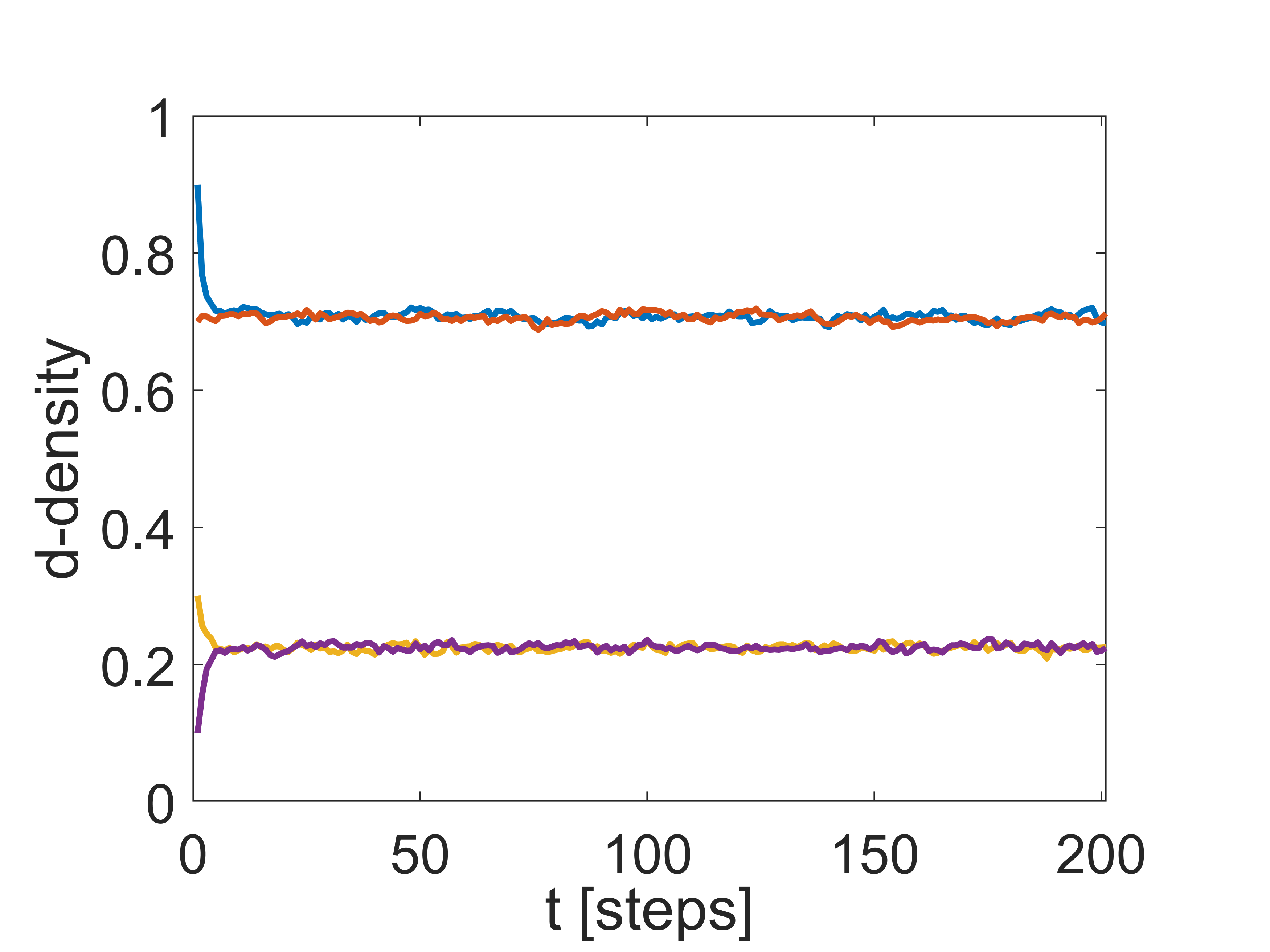}}
    \caption{}
    \label{fig:No_rain_Eff}
    \end{subfigure}
    \begin{subfigure}{.44\textwidth}{\includegraphics[width=1.1\textwidth]{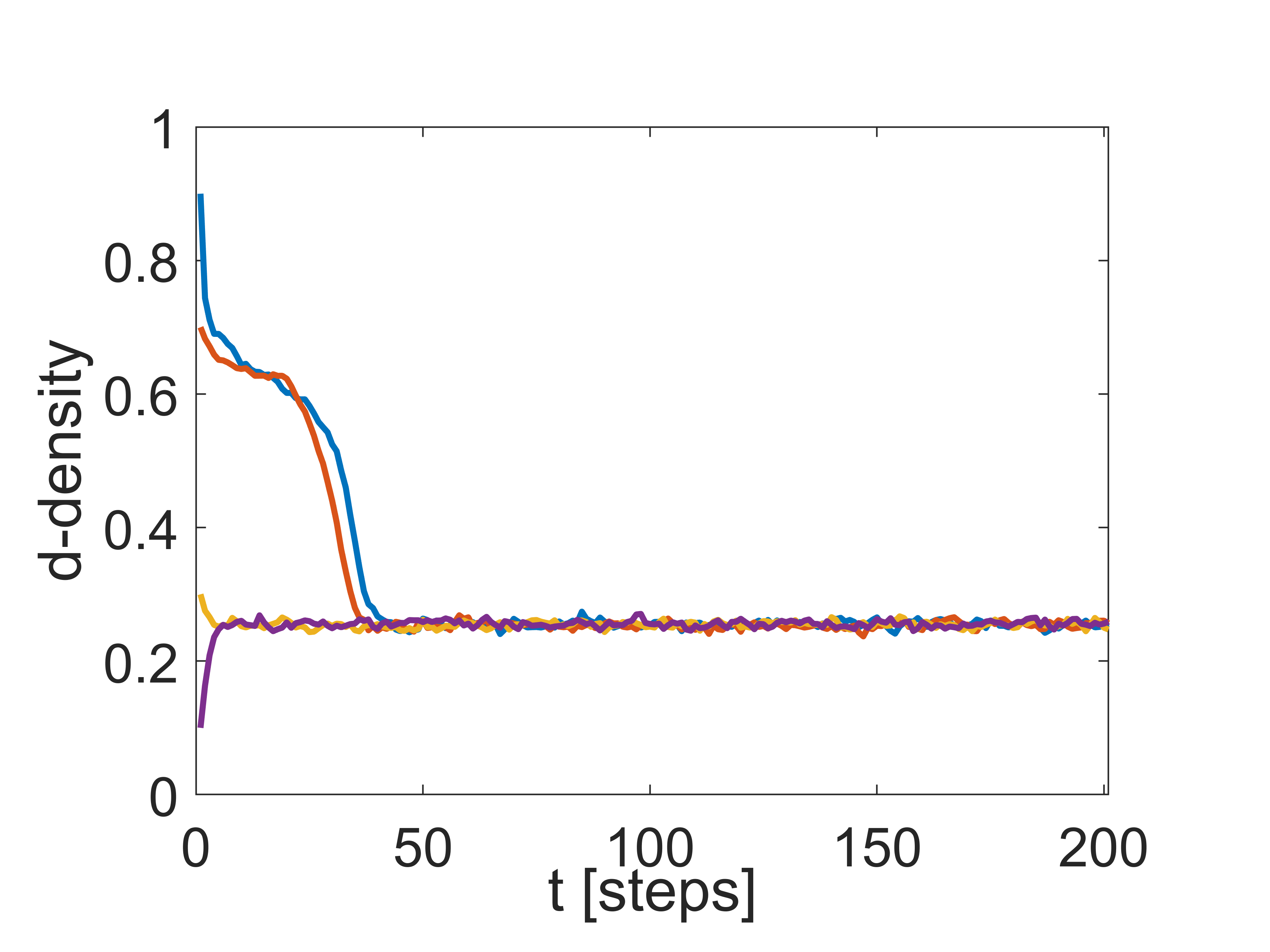}}
    \caption{}
    \label{fig:No_rain_C}
    \end{subfigure}
    \caption{Time-evolution plots of the activity density
for the Erd{\"o}s-R\'enyi type network with $K=10^4$ agents, connectivity probability $p=0.001$, initial activity densities $d_0=0.1, 0.3, 0.7, 0.9$, and activation probabilities $\epsilon = 0.18, 0.2, 0.23, 0.25$. 
For small values of $\epsilon$, two stable states occur (\textbf{a})--(\textbf{c}), while for $\epsilon$ greater than a critical value the network has only one stable state (\textbf{d}).
The high-density solution diminishes (illustrated here for $\epsilon = 0.25$), leading the system, regardless of initial conditions, to converge toward low-density solutions.} 
\label{fig:temp_kozma}
\end{figure}


\section{Characterization of network dynamics using persistent homology}
\label{sec:topology}
As described in the previous section, 
the density of active agents, $d_t$, serves as an order parameter for the state of the network. 
Here, we utilize TDA to create another time series, $r_{\min;t}$, that captures the same long-term behavior of the network. The time series $r_{\min;t}$ is created from a smaller number of agents compared to the activity density time series $d_t$.

\subsection{Order parameter from topological data analysis}
Inspired by the configuration of small world networks \cite{Wat98}, we distribute the agents equidistantly on the unit circle, equipped with the standard round distance function. Then, for a state at a given time $t$, $\mathbf{s}_t=(s_{1;t},s_{2;t},\dots,s_{K;t})$, we consider the Vietoris-Rips complex with the active agents as vertex set and parameter $r$ \cite{Chaz21, Carl09, Ada09}. For small values of the parameter $r$, the first (persistent) Betti number of the associated complex is zero. However, as $r$ increases, the first Betti number eventually becomes one. We define $r_{\min;t}$ as the minimum value of $r$ for which the first Betti number equals one. Note that this minimum value actually depends on the network state at time $t$. Detailed computational procedures are provided in the Supplementary Material. A key finding of this work is that the time series $r_{\min;t}$ captures the same long-term dynamical behavior of the network as the active agent density $d_t$.

\subsection{Witness filtration process}
\label{sec: Witness filtration process}
Computing $r_{\min;t}$ as the time series of the active agent density can be excessively costly when dealing with a large number of agents. To overcome this issue, we use a landmark of agents and apply the \textit{witness filtration process} described in Refs.~\cite{Chaz21, Carl09, Ada09}. This technique aids in extracting statistics on the birth of the first Betti number. For fixed values of the initial activity density $d_0$ and the interaction probability $\epsilon$, the key steps in the witness filtration process can be summarized as follows:

\begin{itemize}
    \item 
    We select an initial microscopic state consistent with $d_0$ at random, and simulate the dynamics. At time $t>0$, we let $Z(t)$ be the set of active agents.
    
    \item From $Z(t)$ we choose a representative landmark set $\mathcal{L}(t)$ using the inductive selection process described in \cite{Javaplex}: First, we fix a suitable size on $\mathcal{L}(t)$, and we choose the first landmark $z_1\in Z(t)$ at random. Then, inductively, given a sequence of $i$ landmarks $z_1,z_2,\dots,z_i$, we choose $z_{i+1}$ as the $z\in Z(t)$ whose distance to $\{z_1,z_2,\dots,z_i\}$ is maximal. 
    (If the choice is not unique, we pick one of the candidates at random.)
    \item For the landmark set $\mathcal{L}(t)$, we perform a filtration process using Lazy Witness simplicial complexes $LW(Z(t), \mathcal{L}(t), r(t))$; see Refs.~\cite{Javaplex, Silva04}. {A visualization of the filtering process is provided in the Supplementary Material.}
    
    \item We record the smallest radius, 
$r_{\min;t}$, where a dimension $1$ hole appears in the filtration. In other words, $r_{\min;t}$ is the minimum value of $r$ for which the first Betti number of the Lazy Witness complex is $1$. Our stochastic agent-based network model generates a unique $r_{\min;t}$ at each time step. 
\end{itemize}
Figure \ref{fig:pip}(A) shows a pipeline of the aforementioned procedure. For a more detailed description of the filtering process, we refer to the Supplementary Material.

We repeat the algorithmic procedure described above multiple times, each time starting with the same initial active agent density $d_0$ and the same interaction probability $\epsilon$, but using different microscopic initial states consistent with $d_0$. For each run, we compute the corresponding $r_{\min;t}$ values and then calculate their average over all runs. The resulting time series 
\begin{equation}
    R_t = \frac{1}{N}\sum_{j=1}^{N} r^{(j)}_{\min;t}. 
    \label{eq: def of R}
\end{equation}
is then compared with the averaged active agent density in the network, 
\begin{equation}
    D_t = \frac{1}{N}\sum_{j=1}^{N} d^{(j)}_t.   
    \label{eq: def of D}
\end{equation}
Here, $N$ denotes the number of network realizations, while $r^{(j)}_{\min;t}$ and $d^{(j)}_t$
represent the minimum radius and the active agent density at time $t$ in the $j$th realization, respectively.

The dynamics of the network concerning the macroscopic density $D_t$ were extensively studied in \cite{Spil10, Spil11} where the authors showed that, depending on the density of the connections and the excitability $\epsilon$, the system exhibits nonlinear behavior, e.g., existence of multiple solutions and critical bifurcation points.
\subsection{Using Betti numbers to express the dynamics of agent-based networks}

In this section, we thoroughly examine the dynamics of activated agents using TDA. Our investigation shows that TDA effectively captures complex network dynamics with significantly fewer agents, demonstrating its distinctive approach and effectiveness.

Indeed, for the construction of Fig.~\ref{fig:Betti_time}, we use the exact states of the network for all time steps identical to Fig.~\ref{fig:temp_kozma}. Then, for each time step we extract the minimal radius $r_{\min;t}$ using the methodology of Section~2.\ref{sec: Witness filtration process}. Figure~\ref{fig:Betti_time} depicts the dynamics of the system with respect to the $r_{\min;t}$ and it resembles an equivalent dynamics as in Fig.~\ref{fig:temp_kozma}. In particular, we observe bistability, i.e., depending on the initial conditions, the system converges either to a high
  or to a low activation state. Additionally, the expected $r_{\min;t}$ deceases its value, 
0.02, 0.017, 0.013 for Fig.~\ref{fig:Betti_time_a}, \ref{fig:Betti_time_b}, \ref{fig:Betti_time_c}), similarly to  Fig.~\ref{fig:temp_kozma} for the low activation state (the expected density is increased).    
\begin{figure}[H]
    \centering
    \begin{subfigure}{.44\textwidth}{\includegraphics[width=1.1\textwidth]{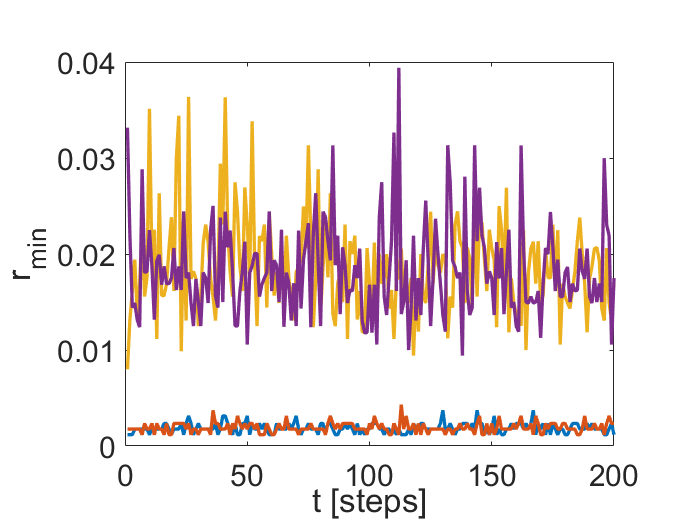}} 
    \caption{}
    \label{fig:Betti_time_a}
    \end{subfigure}
    \begin{subfigure}{.44\textwidth}{\includegraphics[width=1.1\textwidth]{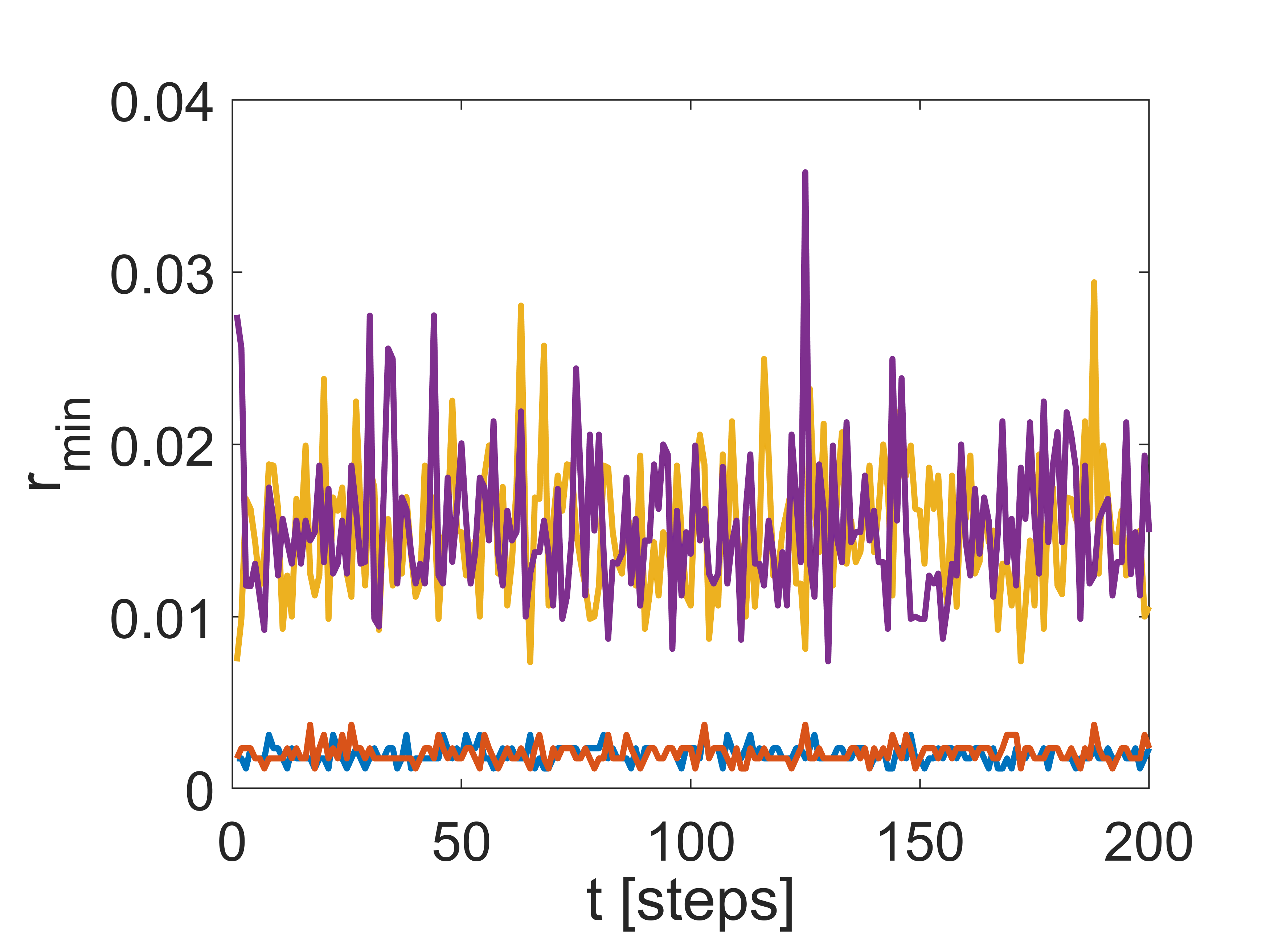}} 
    \caption{}
    \label{fig:Betti_time_b}
    \end{subfigure}
    \begin{subfigure}{.44\textwidth}{\includegraphics[width=1.1\textwidth]{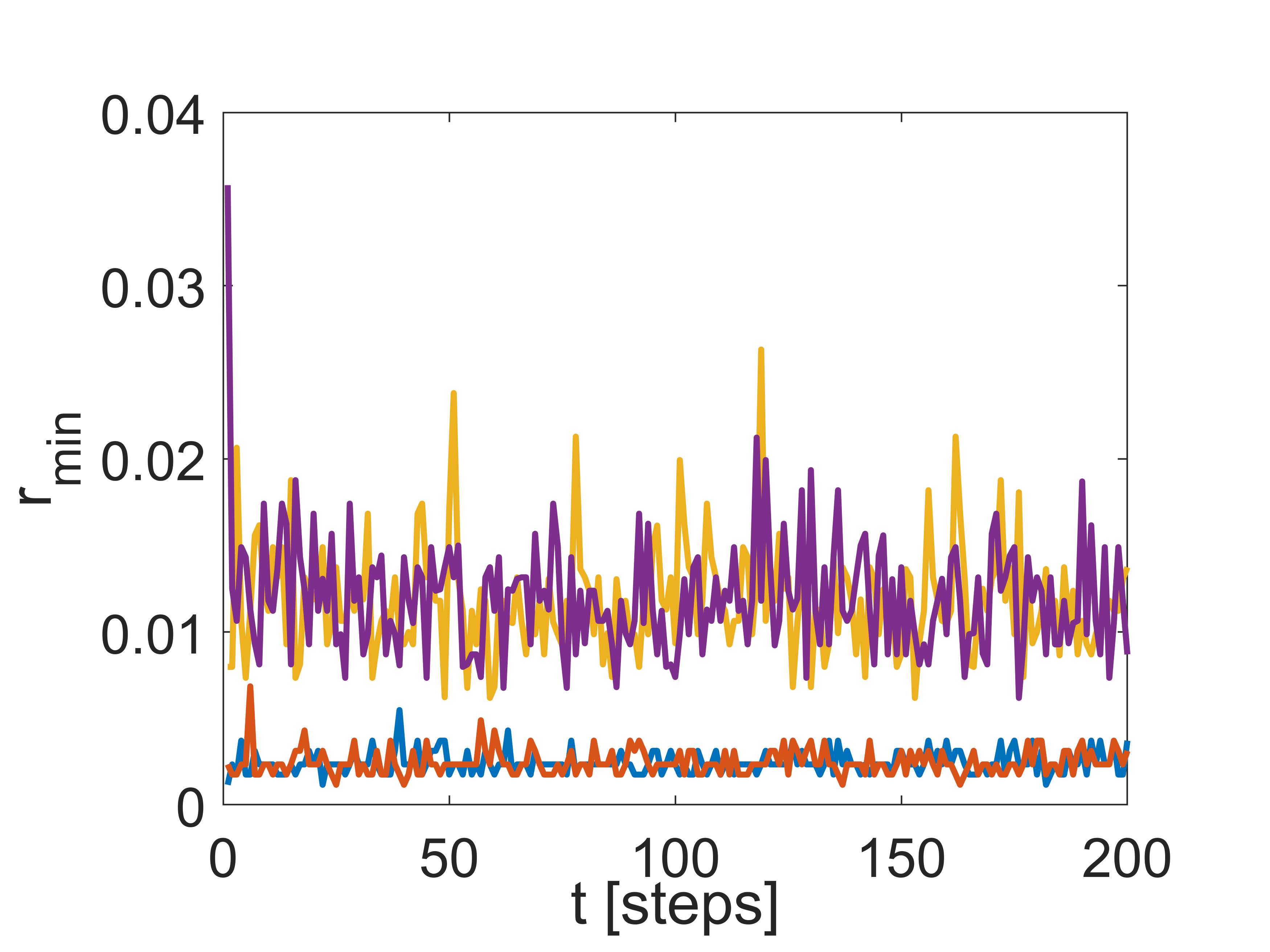}}
    \caption{}
    \label{fig:Betti_time_c}
    \end{subfigure}
    \begin{subfigure}{.44\textwidth}{\includegraphics[width=1.1\textwidth]{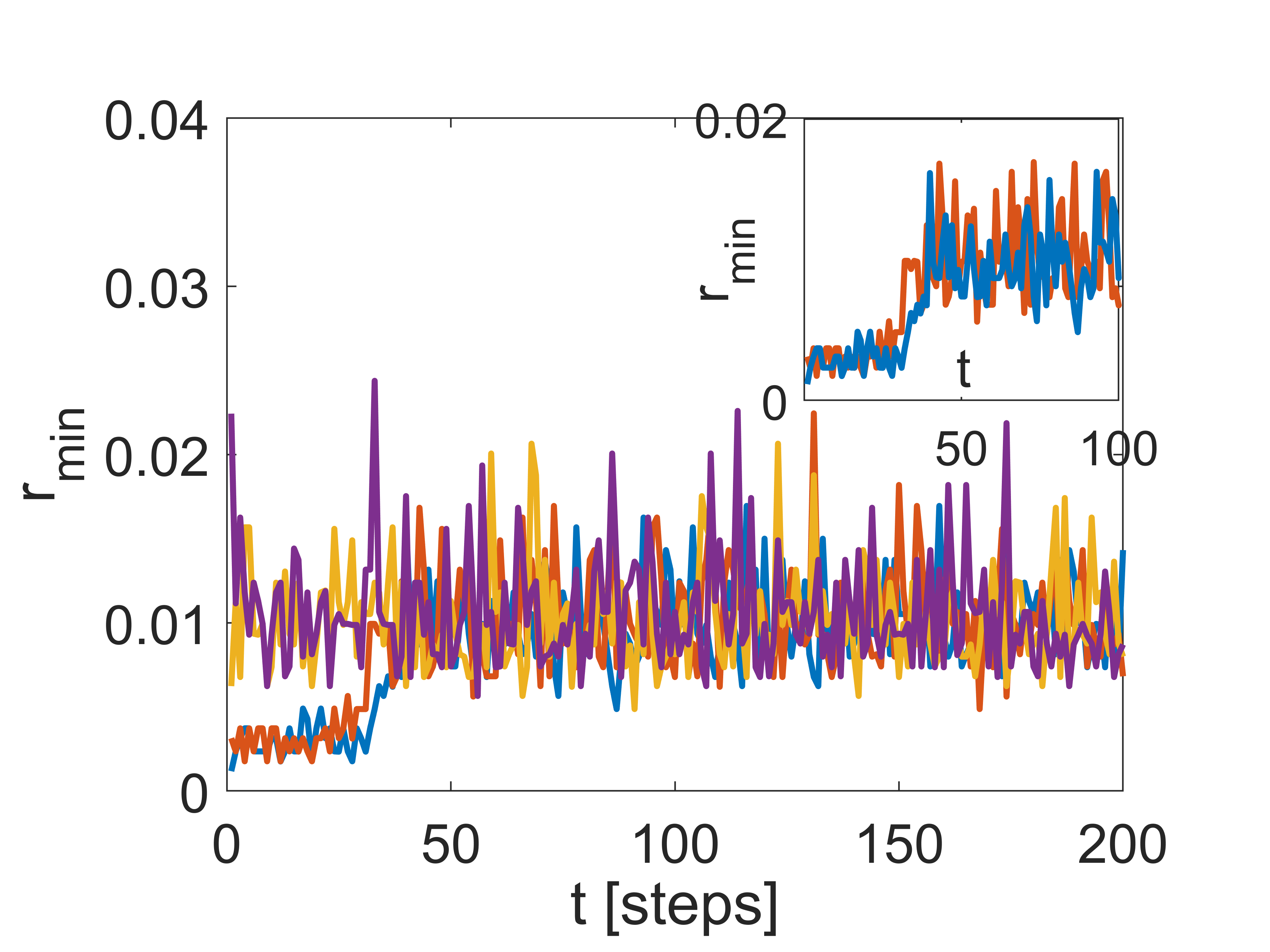}}
    \caption{}
    \label{fig:Betti_time_d}
    \end{subfigure}
    \caption{Using Betti numbers to express the dynamics of agent-based network for $\epsilon=0.18, 0.2, 0.23$ and $\epsilon=0.25$. We plot the minimum radius $r_{\min;t}$ where we obtain $Betti_1=1$. For each case of $\epsilon$ and for each time step, we use the state of the network to compute the $r_{\min;t}$. {Purple and yellow lines represent solutions that are initialized on a low activation state, while red and blue lines are initialized on the high activation state of the network}
\textbf{(a)} $\epsilon=0.18$. The system shows bistability, and depending on the initial conditions, the dynamics converge to the low activation state (the case where $r_{\min;t}$ fluctuates around 0.02) or to the high activation state, where the $r_{\min;t}$ fluctuates around 0.002. 
\textbf{(b)} $\epsilon=0.2$. Similar to (a), the $r_{\min;t}$ in the low activation state fluctuates around 0.017) 
\textbf{(c)} Simulations for $\epsilon=0.23$. The $r_{\min;t}$ in the low activation state fluctuates around 0.013) 
 \textbf{(d)} $\epsilon=0.25$. Now the system, independently of the initial conditions, converges to the low-activation state. The inset shows initialization from high density, which shows a transition to a low-activation state.  
 }
\label{fig:Betti_time}
\end{figure}
\begin{figure}[H]
    \centering
    \begin{subfigure}{.44\textwidth}{\includegraphics[width=1.1\textwidth]{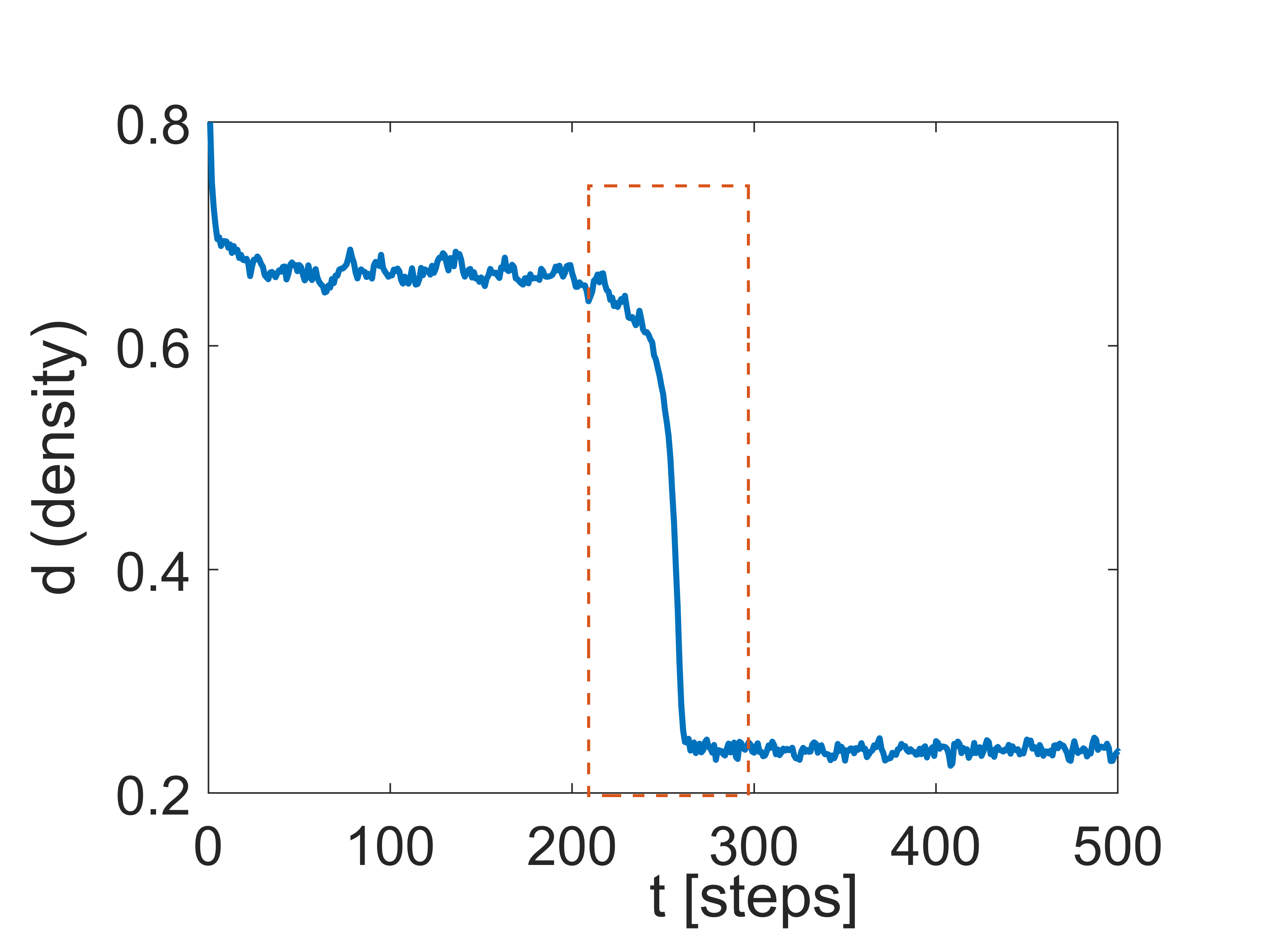}} 
    \caption{}
    \label{fig:Trans_a}
    \end{subfigure}
    \begin{subfigure}{.44\textwidth}{\includegraphics[width=1.1\textwidth]{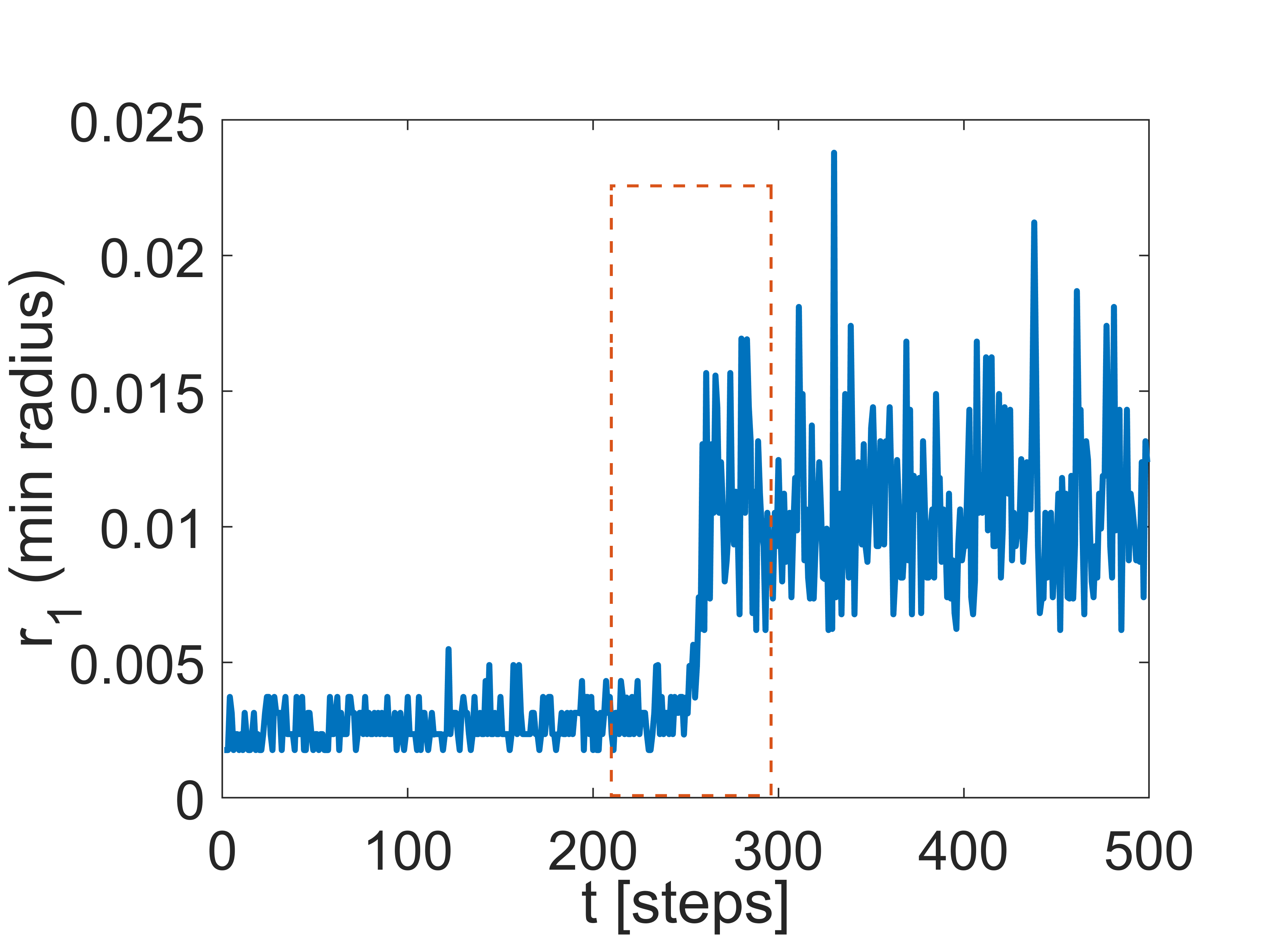}} 
    \caption{}
    \label{fig:Trans_b}
    \end{subfigure}
    \caption{The transient dynamics of the network are observed when simulating the model with $\epsilon=0.239$, initiated from a density of $d_0=0.8$. \textbf{(a)} Illustrates the phase transition in the macroscopic density $d_t$ over time. \textbf{(b)} Notably, a similar transition occurs in the minimal radius $r_{\min;t}$ within the same time horizon.} 
\label{fig:transition} 
\end{figure}



Next, the high-density state, in Fig.~\ref {fig:Betti_time}, results in a low steady state for $r_{\min,t}$. As shown in Fig.~\ref{fig:Betti_time}, the minimum radius $r_{\min,t}$,  for the high activations network state, fluctuates primarily between $0.002$ and $0.004$, i.e., the minimal radius exhibits a notable decrease by approximately one order of magnitude (around $0.002$) compared to the low activation state. 

Finally, in Fig.~\ref{fig:Betti_time_d} we obtain the $r_{\min,t}$ for Fig.~\ref{fig:temp_kozma}d (simulating the system for $\epsilon=0.25$). Independent of the initial conditions, the dynamics with respect to $r_{\min,t}$ converge to a low activation state. Furthermore, the inset of  Fig.~\ref{fig:Betti_time_d} shows that even with a high activation state of the system, the dynamics of $r_{\min,t}$ jumps to the low activation network state. 

As highlighted in Sec.~\ref{sec:model}, depending on the activation probability $\epsilon$, the random network exhibits two stable phases corresponding to high and low-density states. Beyond a critical value of $\epsilon$, only the low-density state persists. The stochastic nature of the majority rule leads to a phase transition from high to low density, as illustrated in Fig.~\ref{fig:Trans_a}. Notably, $r_{\min;t}$ undergoes a similar transition within the same time window, as shown in Fig.~\ref{fig:Trans_b}. Overall, we conclude that both descriptions for density $d_t$ and $r_{\min,t}$ are qualitatively equivalent.
%

\subsection{Macroscopic dynamics}
The macroscopic variables $R_t$ and $D_t$, defined as the averages of $r_{\min;t}$ and $d_t$ (see \eqref{eq: def of R} and \eqref{eq: def of D}), both capture the emergent dynamics of the network's activity. One of the central results of this study is that the time series $R_t$ exhibit stability characteristics that closely mirror those of the macroscopic time series $D_t$ (see Fig.~\ref{fig:betti_macro} and Ref.~\cite{Spil11}). This correspondence highlights the ability of $R_t$ to reflect key dynamical features of the underlying network. Notably, $R_t$ achieves this while requiring significantly less information than $D_t$. {Precisely, we are using a set of $50$ landmark agents to compute $R_t,$ while in contrast, for calculating $D_t$ we need $10000$ agents.}

Figure~\ref{fig:betti_macro} illustrates the dynamics of $R_t$ and $D_t$ for both low and high initial activity densities, respectively. When the activation probability is below the critical value, the network settles into two stable states, one corresponding to a high initial activity density and the other to a low initial activity density, as shown in Figs.~\ref{fig:betti_macro_low} and~\ref{fig:betti_macro_low_time}. However, when the activation probability exceeds the critical value, the system evolves toward a single steady state regardless of its initial activity density, as illustrated in Figs.~\ref{fig:betti_macro_high} and~\ref{fig:betti_macro_high_time}.



\section{Exploring the network dynamics using the Equation-Free Method}
\label{sec:results}

Our next objective is to formulate an evolutionary equation for the macroscopic time series $R_t$ of the form
\begin{equation}
R_{t+T} = F_{T}(R_t;\epsilon),
\label{eq: Time-stepper}
\end{equation}
where $T$ denotes a macroscopic time step. Such a formulation allows for system-level analysis, including the identification of steady states, examination of stability and control properties, and the study of bifurcations. However, obtaining an explicit closed-form expression for the time-stepper $F_{T}$ is generally a highly nontrivial task. To address this, we employ an algorithmic framework known as the Equation-Free Method (EFM), which provides a numerical approximation of the time-stepper; see Ref.~\cite{Gear03}.

\begin{figure}[H]
    \centering
    \begin{subfigure}{.44\textwidth}{\includegraphics[width=1\textwidth]{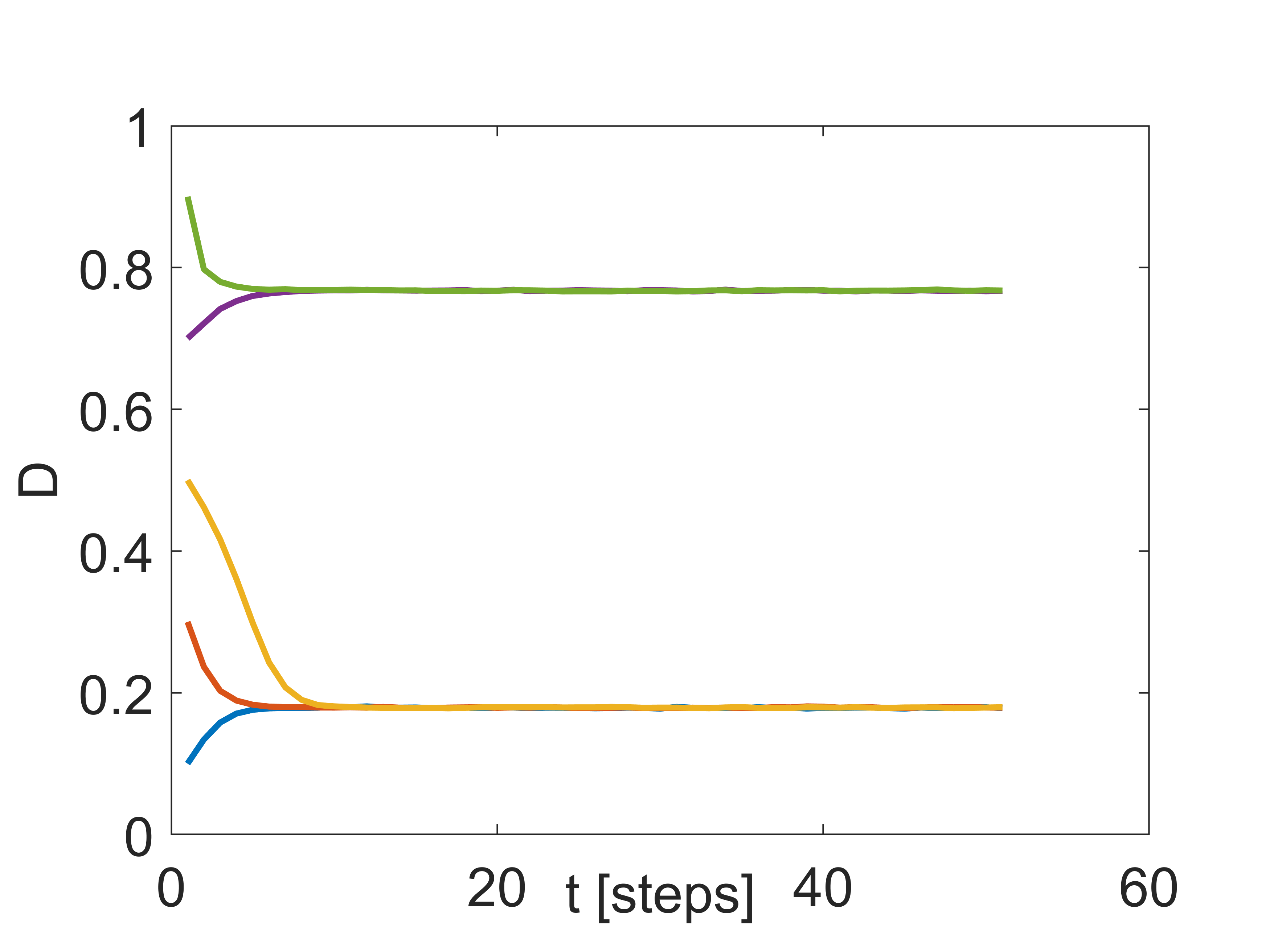}} 
    \caption{The evolution of $D_t$ for $\epsilon = 0.2$.}
    \label{fig:betti_macro_low}
    \end{subfigure}
    \begin{subfigure}{.44\textwidth}{\includegraphics[width=1\textwidth]{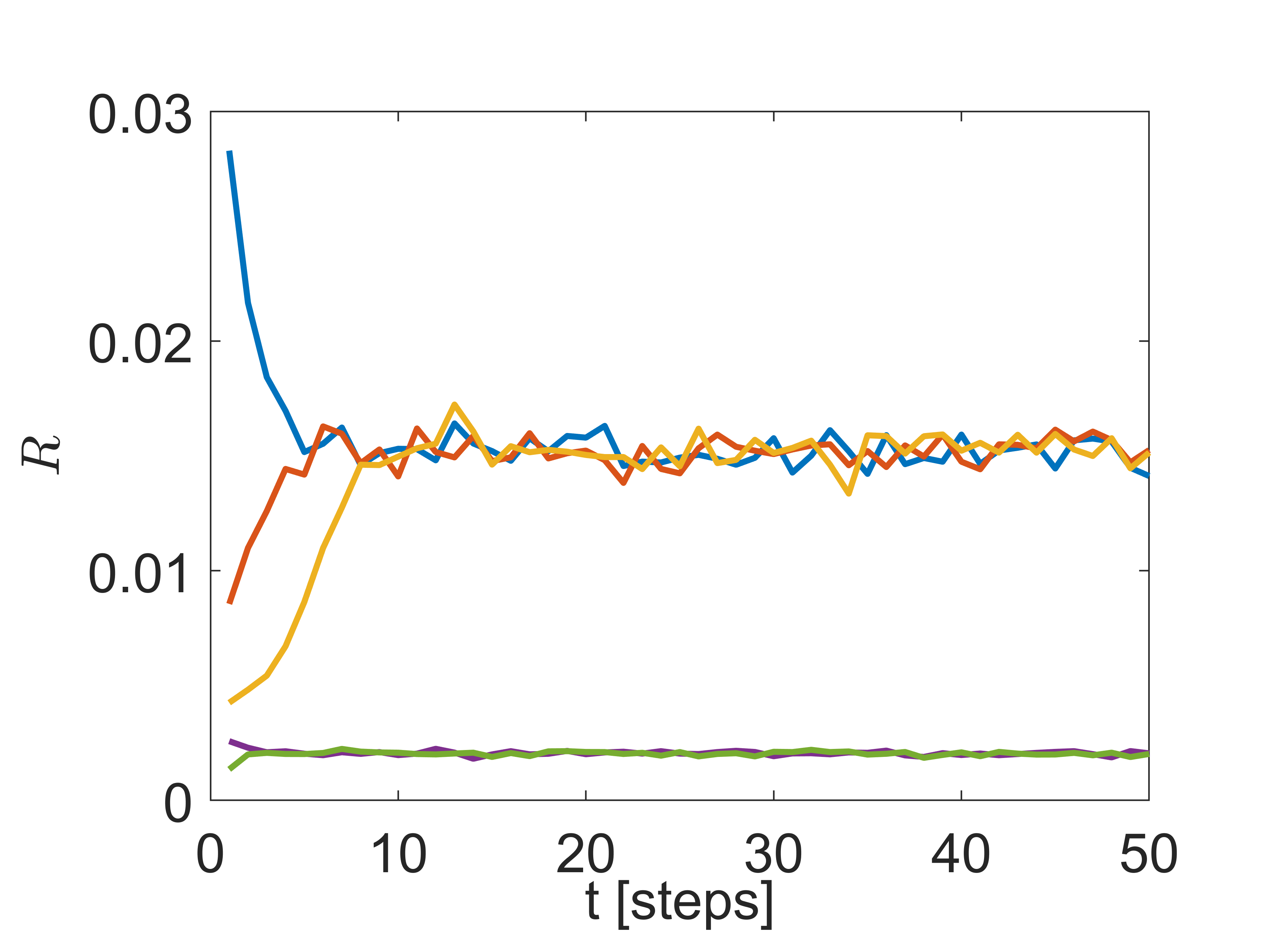}} 
    \caption{The evolution of $R_t$  for $\epsilon = 0.2$.}
    \label{fig:betti_macro_low_time}
    \end{subfigure}
    \begin{subfigure}{.44\textwidth}{\includegraphics[width=1\textwidth]{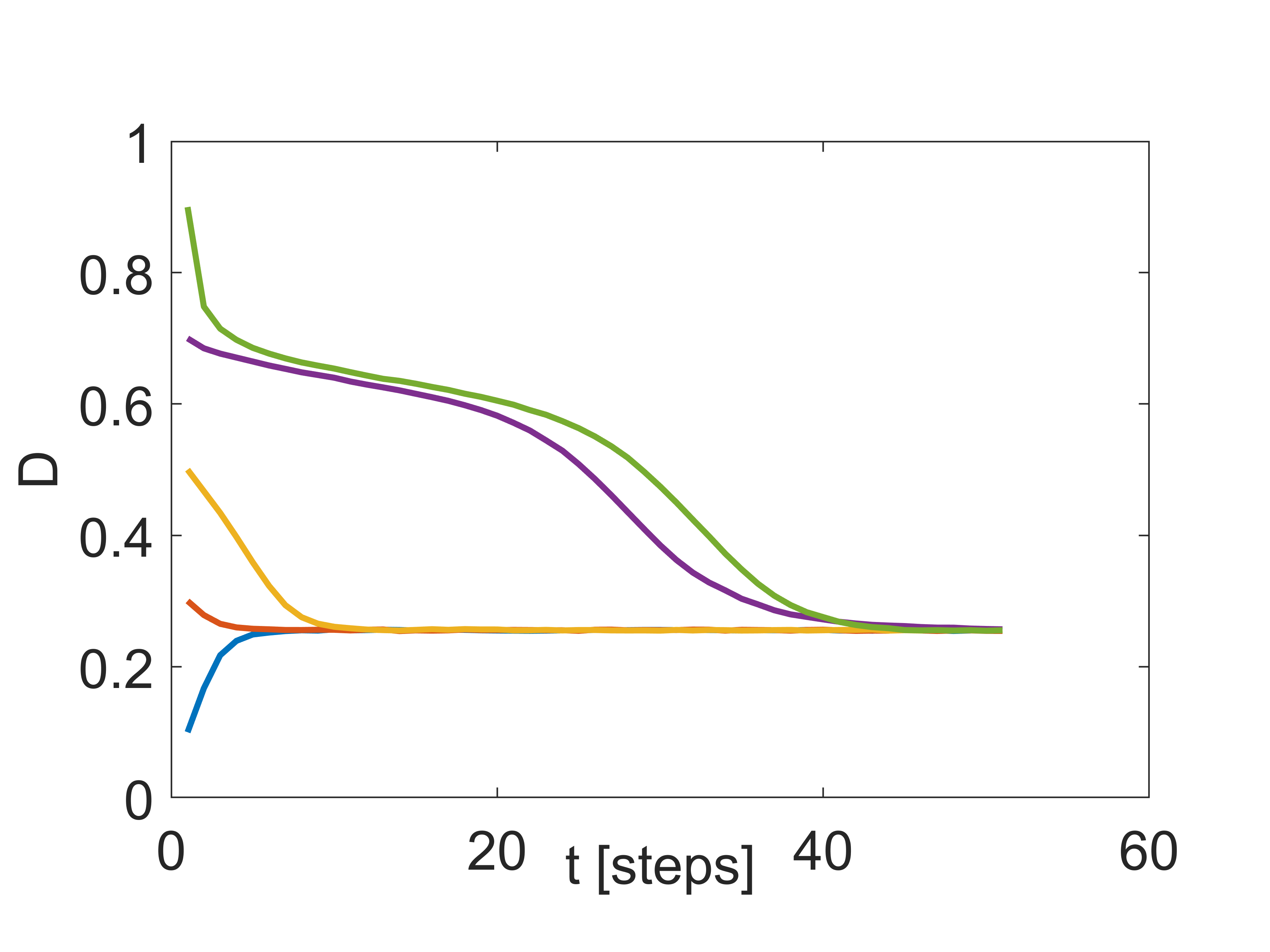}}
    \caption{The evolution of $D_t$  for $\epsilon = 0.25$.}
    \label{fig:betti_macro_high}
    \end{subfigure}
    \begin{subfigure}{.44\textwidth}{\includegraphics[width=1\textwidth]{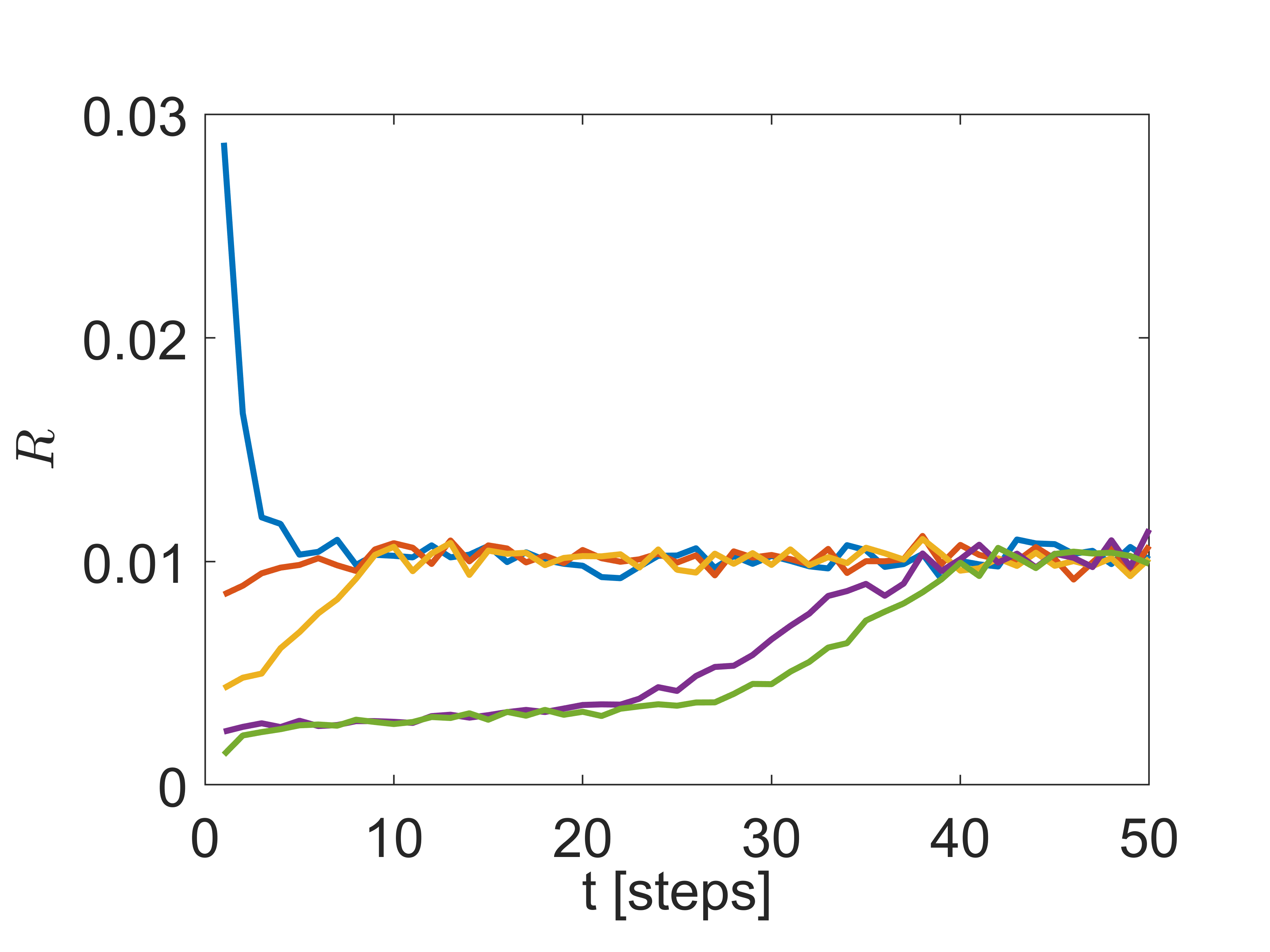}}
    \caption{The evolution of $R_t$ for $\epsilon = 0.25$.}
    \label{fig:betti_macro_high_time}
    \end{subfigure}
    \caption{The evolution of macroscopic variables $D_t$ and $R_t$ is shown for different initial activation densities and activation probabilities. For $\epsilon = 0.2 < \epsilon_c$, the system settles into two stable states: a low-activity state with $D = 0.18$ and $R = 0.017$, and a high-activity state with $D = 0.78$ and $R = 0.002$ (see (a) and (b)). In contrast, for $\epsilon = 0.25 > \epsilon_c$, the system evolves toward a single steady state regardless of its initial activity density, characterized by $D = 0.25$ and $R = 0.011$ (see (c) and (d)).}
\label{fig:betti_macro}
\end{figure}


Traditionally, generating a macroscopic time series requires multiple independent simulations, each initialized with different microscopic states consistent with the same macroscopic condition. The resulting trajectories are averaged to estimate the macroscopic dynamics, but this process is computationally expensive. The EFM addresses this by constructing the entire macroscopic time series in a single execution. This not only improves efficiency but also yields an explicit dynamical equation for macroscopic evolution, enabling rigorous stability and bifurcation analysis of the underlying network dynamics.


The EFM involves three operators, the lifting operator $L$, the evolution operator $\Phi_{T}$, and the restriction operator $M$. These operators do the following:
\begin{itemize}
    \item The lifting operator $L: \mathbb{R}^+ \to \{0, 1\}^K$ maps a given positive value $R$ to a randomly generated network configuration $\bold{s}$ such that the minimum radius $r_{\min}$ computed from $\bold{s}$ is equal to $R$; see Secs.~3.\ref{sec:lifting_operator_1} and 3.\ref{sec:lifting_operator_2}.
    \item The evolution operator $\Phi_{T} : \{0,1\}^K\times (0,0.5) \to \{0,1\}^K$ maps a network configuration $\bold{s}$ to a new configuration $\Phi_{T}(\bold{s};\epsilon)$ by evolving it over a time interval of length $T$ according to the majority rule with activation probability $\epsilon$, as detailed in Sec.~1.\ref{sec: Description of the network}.
    
    \item The restriction operator $M:\{0,1\}^K\to \mathbb{R}^+$ takes a network configuration and calculates the corresponding $r_{\min}$.
\end{itemize}
\begin{figure}[H]
\centering
\hspace{-0.5cm}%
\includegraphics[width=0.6\linewidth]{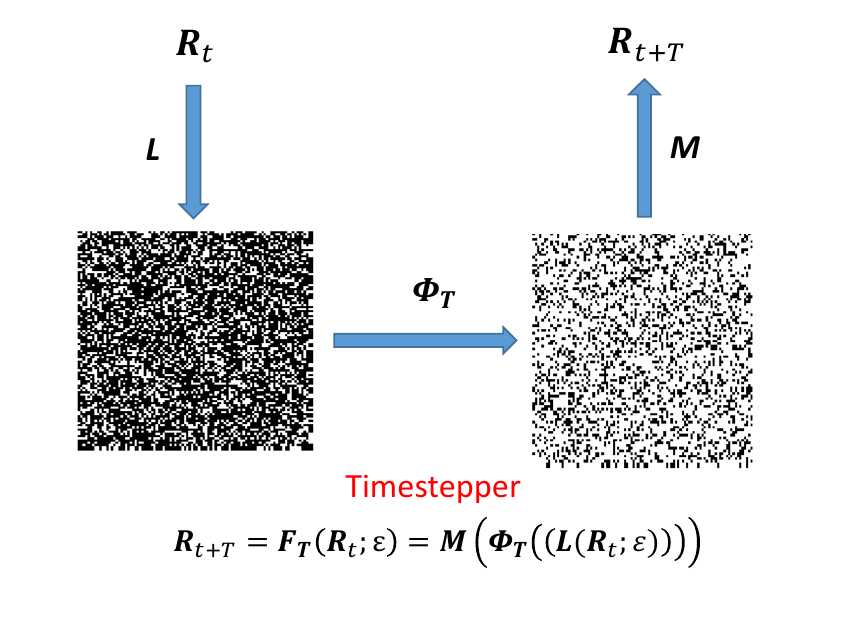}
\caption{A schematic representation illustrating the Equation-Free approach is presented. In this method, we construct the time stepper, which serves as a numerical discretization of the macroscopic dynamics. The time stepper is comprised of the lifting operator $L$, responsible for evolving the microscopic rules, and the restriction operator $M$.  In our context, the macroscopic variable is denoted as $R_t$, representing the minimal radius in the filtration process associated with the emergence of circles.} 
\label{fig:RHS}           
\end{figure}
We now define an evolution law (a time-stepper) as follows,
\begin{equation}
   \tilde R_{(n+1)T}=M(\Phi_T(L(\tilde R_{nT});\epsilon)),\quad \tilde{R}_0=R_0,
   \label{eq:dyn_sys_neuron_1}
\end{equation}
where $n=0,1,2,\dots$. 
\begin{figure}[H]
    \centering
    \begin{subfigure}{.44\textwidth}{\includegraphics[width=1.0\textwidth]{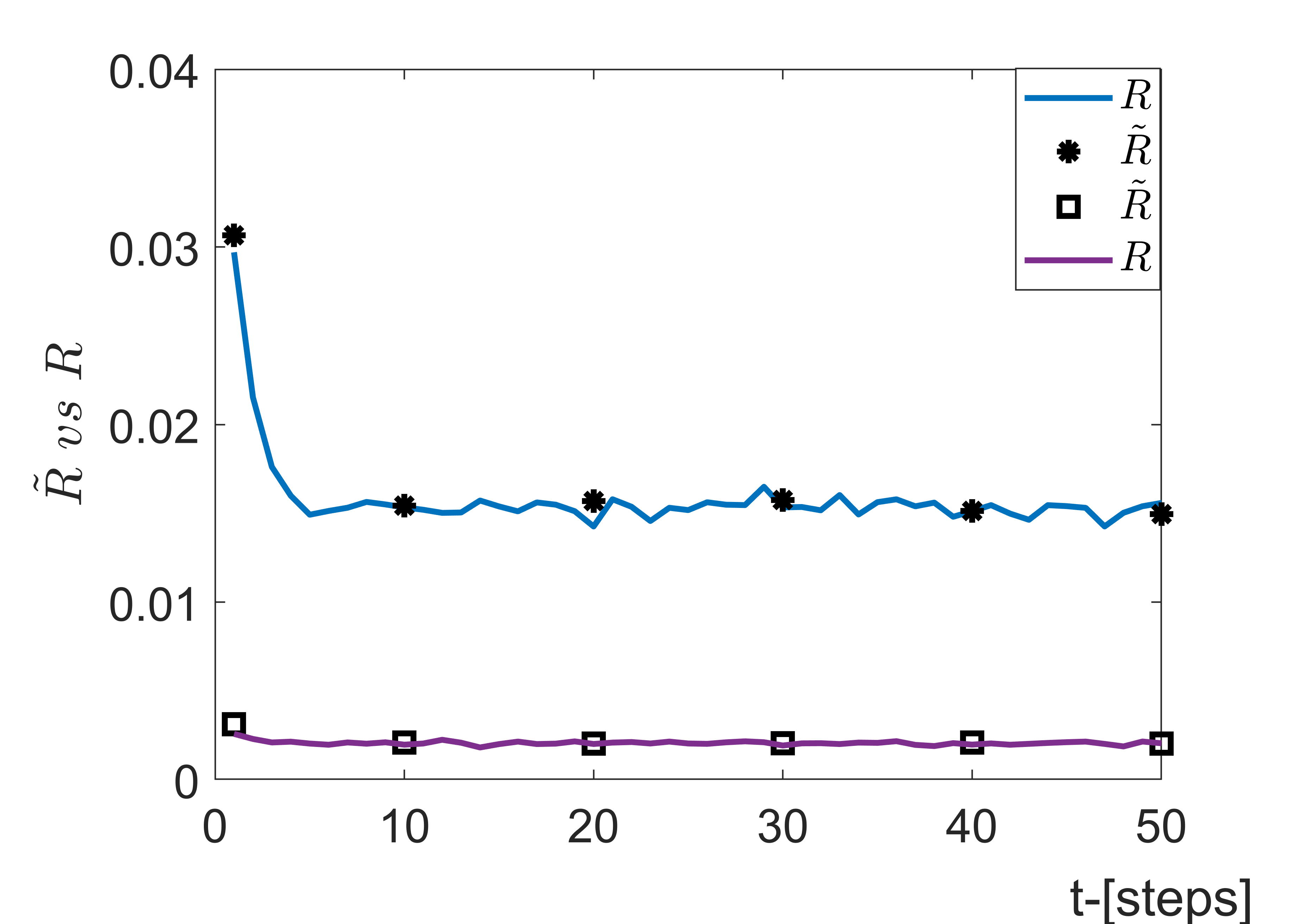}} 
    \caption{}
    \label{fig:rtilda_time}
    \end{subfigure}
    \begin{subfigure}{.44\textwidth}{\includegraphics[width=0.97\textwidth]{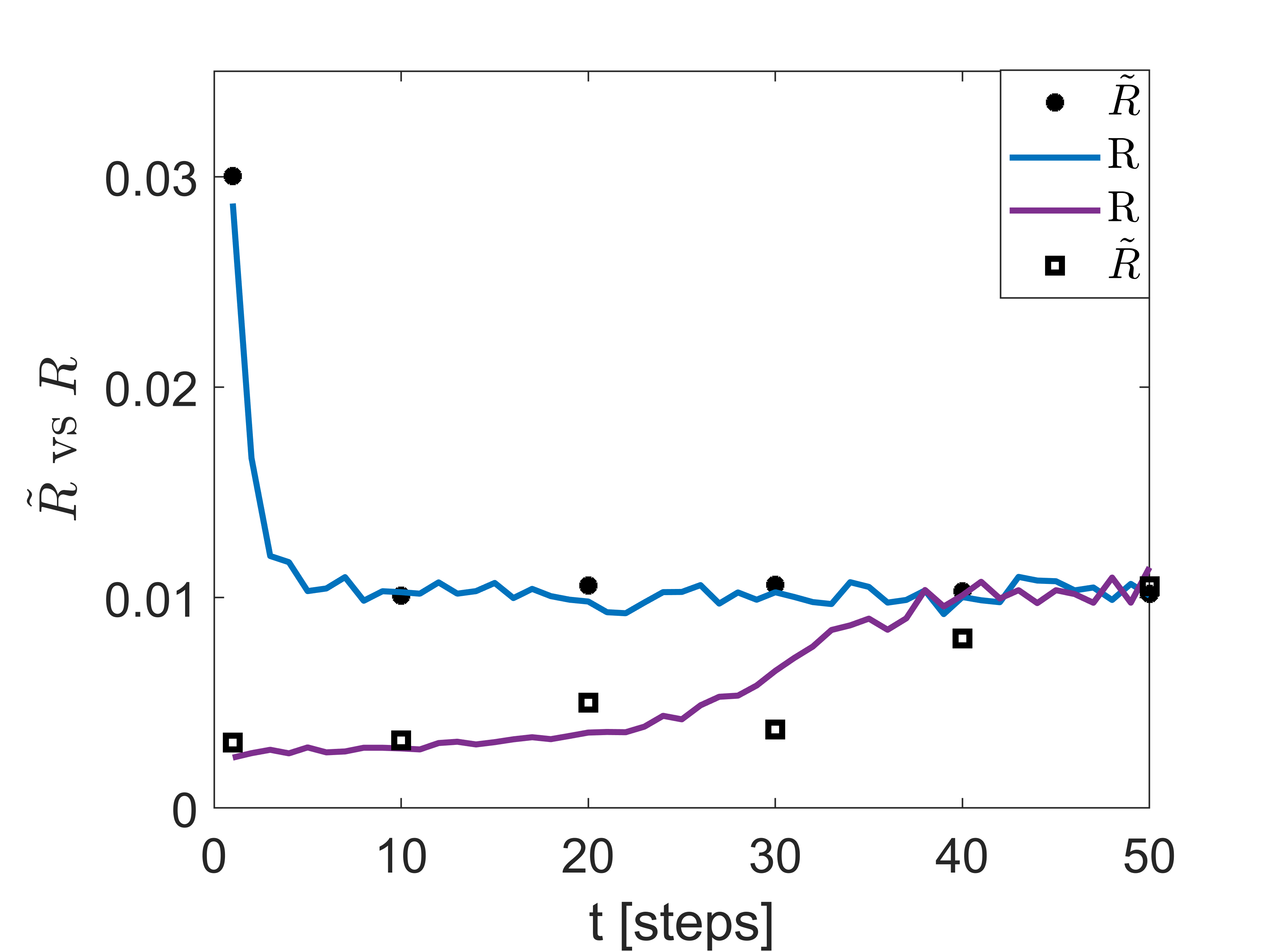} }
    \caption{}
    \label{fig:manifold_timestepper}
    \end{subfigure}
    \caption{The diagrams illustrate the dynamics of the coarse time-stepper $\tilde{R}_t$, defined in~\eqref{eq:dyn_sys_neuron_1}, alongside the time series $R_t$. Black dots indicate values of $\tilde{R}_t$ at discrete times $t = nT$ with $T = 10$, while the continuous curves represent $R_t$. Panel (a) corresponds to the bistable regime with $\epsilon = 0.2$ and should be compared to Figure~\ref{fig:betti_macro_low_time}. Panel (b) shows the monostable case for $\epsilon = 0.25$, corresponding to Figure~\ref{fig:betti_macro_high_time}. The coarse time series $\tilde{R}_t$ was generated using the geometric information for the lifting procedure described in Section~3.\ref{sec:lifting_operator_1}.}
\label{fig:rtilda_and_r}
\end{figure}
Under the main assumptions of EFM \cite{Gear03, Kev09, Sieb18}, the time series $\tilde{R}_{t}$ closely approximate $R_t$ on the macroscopic time scale $T$ (see Fig.~\ref{fig:rtilda_and_r}). Therefore, from this point forward, we will treat the two interchangeably.

A representative bifurcation diagram for $R_t$ generated using the EFM is presented in Fig.~\ref{fig:RHS1}, offering further insight into the qualitative transitions observed in the system's topological evolution. {We compute steady state solutions using the Newton-Raphson scheme and derivative information for stability analysis, see also the caption of Fig.~\ref{fig:RHS1}}.

\subsection{The lifting operator using geometric information}
\label{sec:lifting_operator_1}
In general, constructing the lifting operator (typically a one-to-many correspondence problem) is one of the most challenging aspects of the EFM, often necessitating an optimization procedure \cite{SIET11, Sieb17, Sieb18, PATSA23}. In our setting, however, TDA offers a natural and canonical framework for defining $L$, which we now describe.

Consider the distribution of agents on the unit circle, as detailed in Sec.~2.\ref{sec: Witness filtration process}. 
For a given input $R>0$, let $\delta(R)$ represent the distance closest to $2R$ among all the pairwise distances between agents. 
Thus, $\delta(R)/2\approx R$. To construct the subset of agents $Z=\{z_1,z_2,\dots,z_Q\}\subset\{a_1,a_2,\dots,a_K\}$, 
proceed as follows: start with $z_1=a_1$, then iterate through the remaining agents, adding an agent to 
$Z$ if its distance to the last added agent is $\delta(R)$. The process is illustrated in Fig.~\ref{fig:geo_lift}

The network state $L(R)=\bold{s} =(s_1,s_2,\dots,s_K)$ is then defined as follows: $s_j=1$ if and only if agent $a_j$ belongs to $Z$, 
and $s_j=0$ otherwise. 
If we compute $r_{\min}$ for the state $L(R)$, we find that
\begin{equation}
    r_{\min}(L(R))=\delta(R)/2\approx R.
\end{equation}

\begin{figure}[H]
    \begin{center}
    \includegraphics [width=0.5\linewidth]{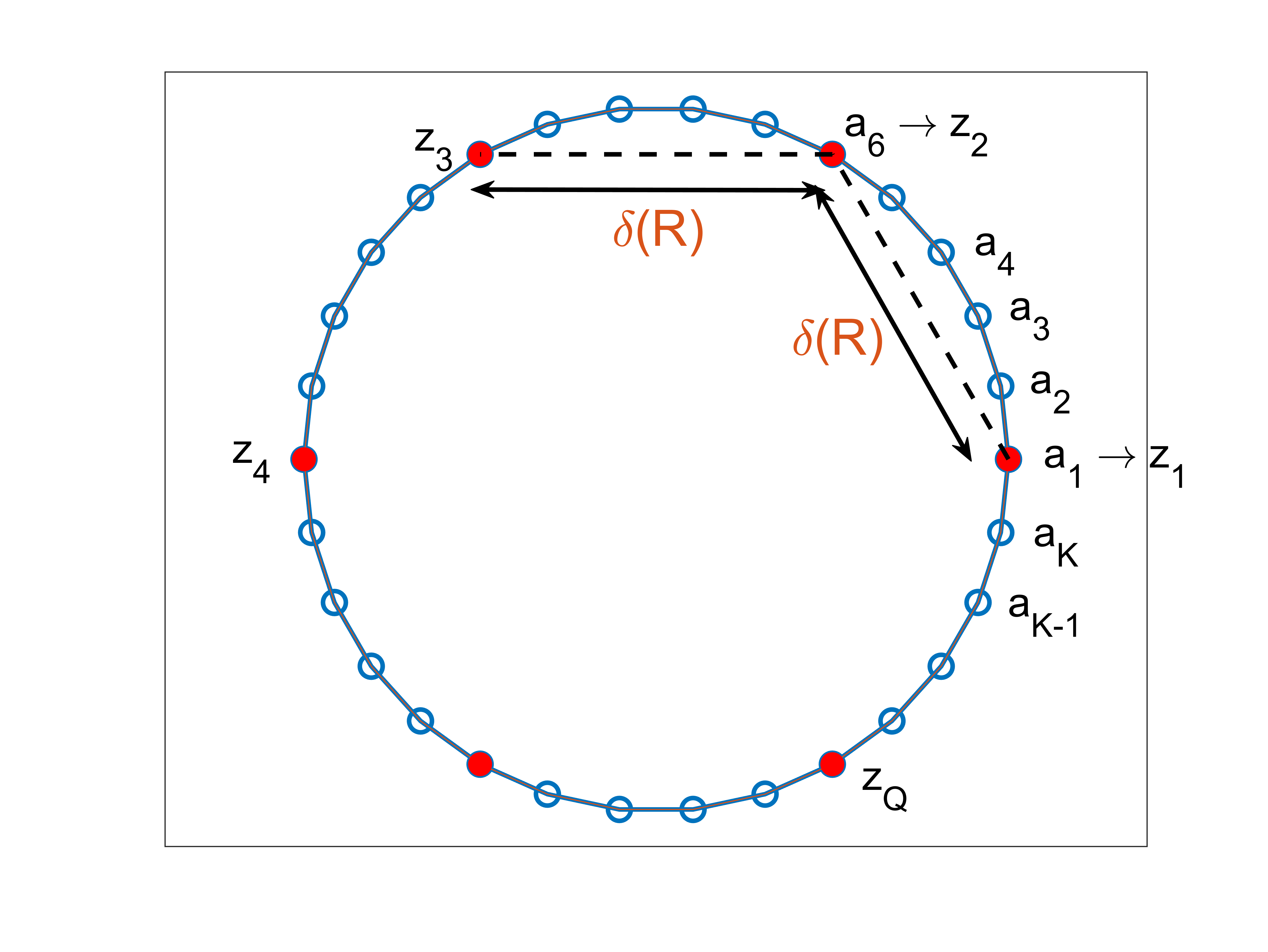}
    \end{center}
    \caption{Construction of the lifting operator using geometric information derived from TDA, as detailed in Section~3.\ref{sec:lifting_operator_1}.}
\label{fig:geo_lift}
\end{figure}

\subsection{The lifting operator using simulated annealing}
\label{sec:lifting_operator_2}
A more general (though not necessarily more efficient) approach for constructing the lifting operator involves the use of simulated annealing (see Refs.~\cite{Kir83, Spil10}). Simulated annealing is a stochastic optimization technique that minimizes the energy function
$E=|R-R_{\text{target}}|$, where $R_{\text{target}}$ denotes the desired (target) value of the mean radius, and $R$ is the mean radius corresponding to the current network state.

At each iteration, a candidate microscopic network configuration is generated, leading to a new value of $R$. This is achieved by evolving the system under the majority rule dynamics for a short time horizon 
starting from a random initialization. The updated configuration yields a new mean radius, denoted $R_1$, and the associated energy is computed as $E_1 = |R_1 - R_{\text{target}}|$.

The new configuration is then either accepted or rejected based on the Metropolis criterion (see Algorithm~\ref{alg:algorithm1}):
\begin{itemize}
\item If $E_1 < E_0$, where $E_0$ is the energy of the previous state, the new configuration is accepted unconditionally.
\item If $E_1 > E_0$, the configuration is accepted with probability $e^{-(E_1 - E_0)/\theta_n}$; otherwise, it is rejected.
\end{itemize}

If a configuration is rejected, a new initialization is sampled, and the procedure is repeated until a configuration is accepted. Once a state is accepted, it is retained as the new microscopic realization. This iterative process continues until a termination criterion is satisfied, for instance, when the energy falls below a pre-defined threshold, such as $E=|R-R_{\text{target}}|<0.002$.

\begin{algorithm}[H]
\caption{Lifting operator: transforming $R$ to a network state $\bold{s}$.}
\begin{algorithmic}[1]
\State{\textbf{Input:} $R_{\text{target}}$ value of minimal filtration radius.}
\State {Assign the initial pseudo-temperature, denoted as $\theta_0$, and choose the annealing scheme for the pseudo-temperature reduction.}
\State {\textbf{Do until convergence}}

\begin{itemize}
  \item Start with an initial microscopic network configuration $\mathbf{s}_0.$
 \item Compute the minimum filtration radius at which $Betti_1=1$, denoted as $R_0$, and the corresponding error function $E_0 = |R_0 - R_{\text{target}}|.$
\item Evolve the microscopic configuration of the network using the majority rule for three steps. Calculate the updated error function $E_1 = |R_1 - R_{\text{target}}|.$
    \item {Apply the Metropolis criterion to either accept or reject the new configuration $\mathbf{s}_{1}$:} 
    \begin{itemize}
    \item If $E_{1}<E_{0}$, accept the new configuration.
    \item If $E_1 >E_0$, accept the new configuration with probability $e^{-(E_1-E_0)/ \theta_n}$; otherwise reject it.
    \end{itemize}
    
   \item {Decrease the system's pseudo-temperature following the chosen annealing schedule. In this instance, opt to reduce the temperature after $200$ iterations using the formula $\theta_{n+1}=0.9~\theta_{n}$. }
\end{itemize} 
\State {\textbf{End Do}}
\end{algorithmic}
\label{alg:algorithm1}
\end{algorithm}
\subsection{Bifurcation and stability analysis}
As described above, repeated application of the operators $L$, $\Phi_{T}$, and $M$, 
facilitates a numerical approximation
of the time-stepper $F_{T}$ in \eqref{eq: Time-stepper}.

Figure~\ref{fig:Identific} shows the graphs of the function $G_{T}(R;\epsilon) = F_{T}(R;\epsilon) - R$ for $R = R_t$ and $\epsilon = 0.16, 0.19, 0.22, 0.24$. Notably, the overall shape of $G_{T}$ remains qualitatively consistent across all values of $\epsilon$. These characteristic patterns provide meaningful insights into the macroscopic nonlinear dynamics of the network, illustrating how the evolution operator adapts to different activation probabilities.

In particular, the points where the graph of $G_{T}(R;\epsilon)$ intersects the $R$-axis correspond to steady-state solutions. To locate these fixed points, we employ a Newton--Raphson solver applied to $G_{T}(R;\epsilon)$. 

For values of $\epsilon$ below the critical threshold, the system exhibits three steady-state solutions, as illustrated in Fig.~\ref{fig:Identific}. A stability analysis based on the sign of the derivative of $G_{T}(R;\epsilon)$ at each fixed point reveals that the leftmost and rightmost solutions, corresponding to low and high activation states, respectively, are stable, while the intermediate fixed point is unstable.

As $\epsilon$ increases past the critical value, the high-activation stable state and the intermediate unstable state coalesce and annihilate in a saddle-node bifurcation. Beyond this point, the system admits a single stable steady state associated with low activity. This bifurcation scenario is captured in the diagram shown in Fig.~\ref{fig:biftot}.

\begin{figure}[htb]
    \centering
    \begin{subfigure}{.44\textwidth}{\includegraphics[width=1.1\textwidth]{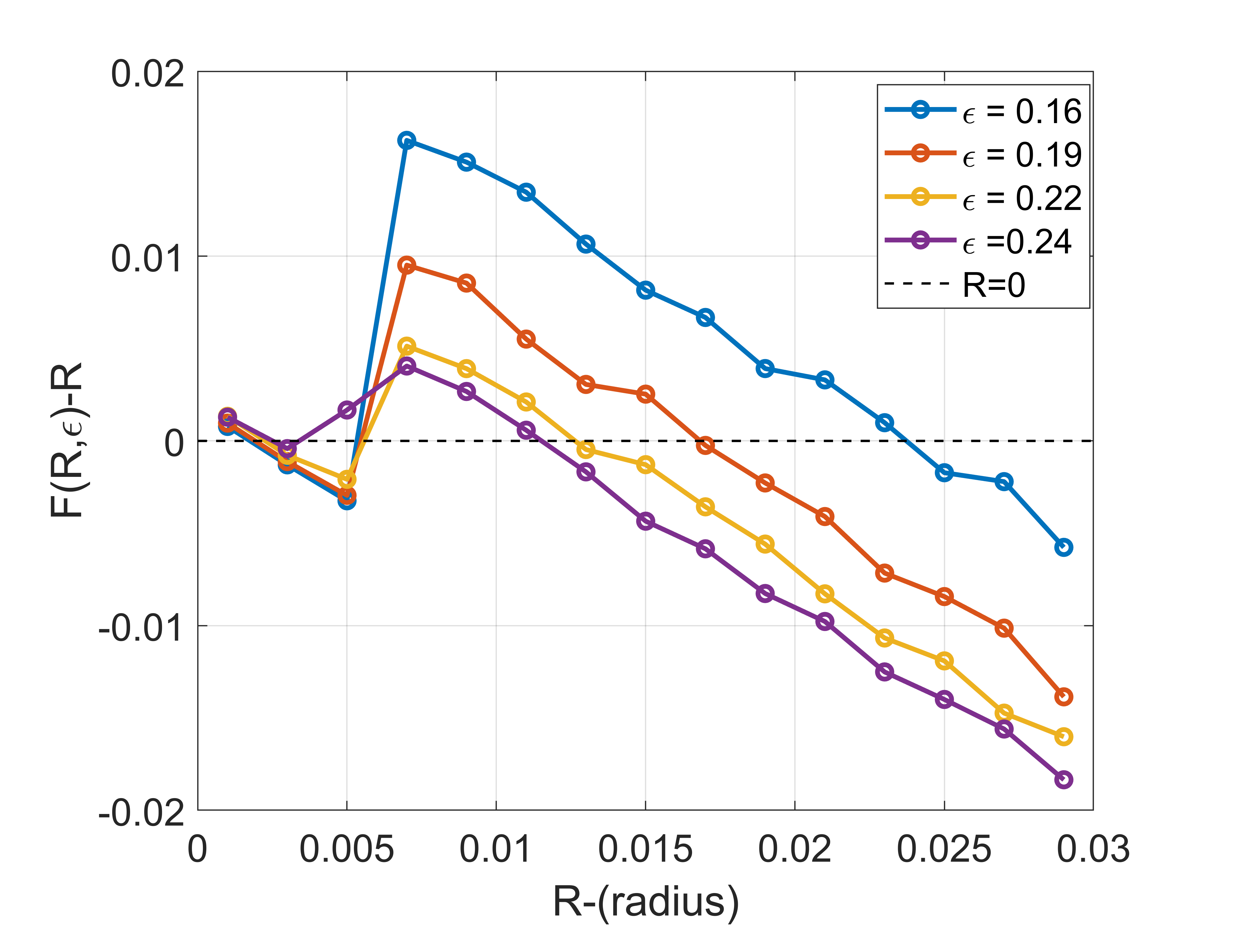}} 
    \caption{}
    \label{fig:Identific}
    \end{subfigure}
    \begin{subfigure}{.41\textwidth}{\includegraphics[width=1.1\textwidth]{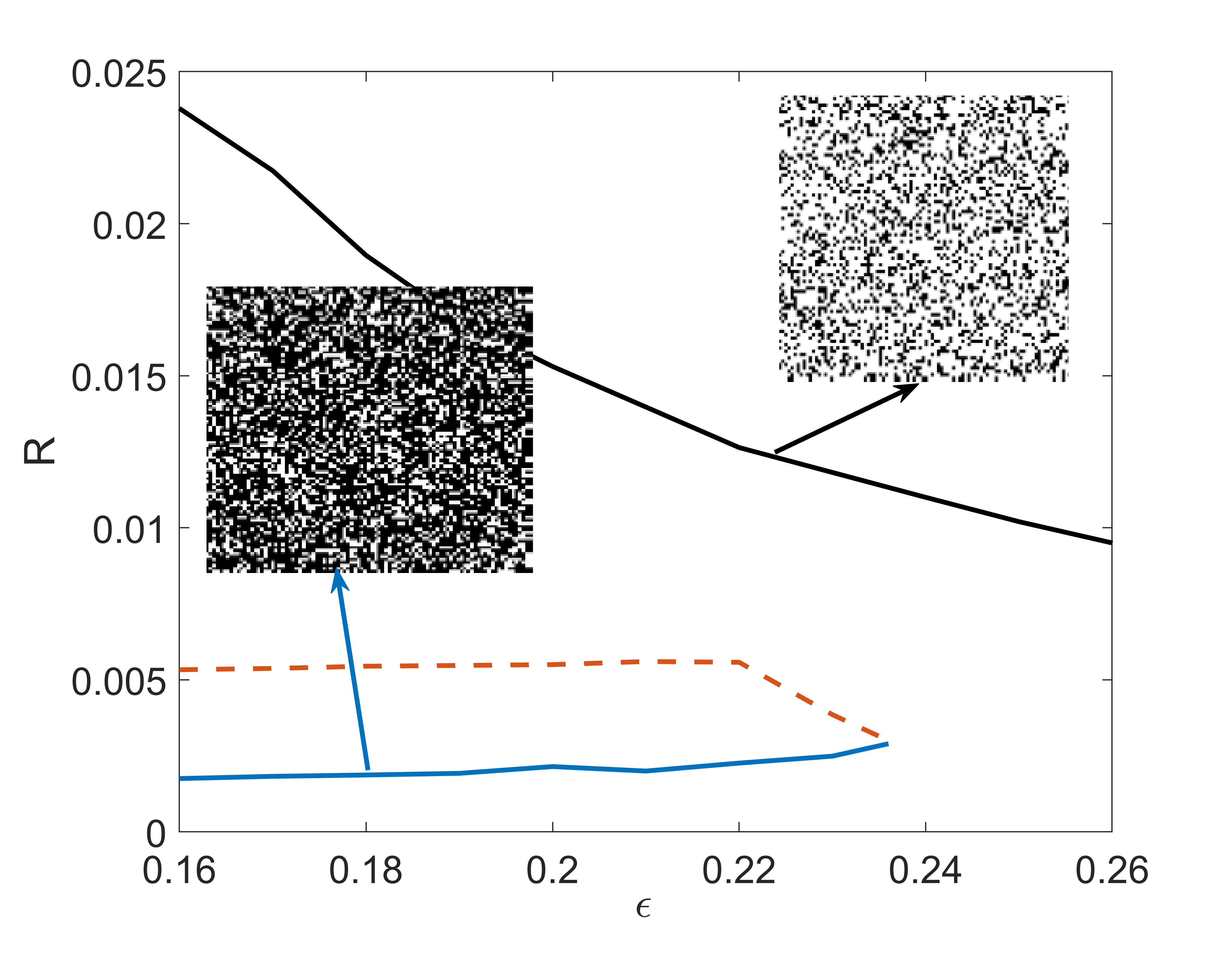} }
    \caption{}
    \label{fig:biftot}
    \end{subfigure}
    \caption{System identification using EFM and TDA.
(\textbf{a}) Plot of the function $G_T(R;\epsilon) = F_T(R;\epsilon) - R$ for $\epsilon = 0.16, 0.19, 0.22$ (below the critical value $\epsilon_c$), and $\epsilon = 0.24, 0.26$ (above $\epsilon_c$). For each $\epsilon < \epsilon_c$, the system exhibits three steady-state solutions at $R_1 \approx 0.002$, $R_2 \approx 0.005$, and $R_3 > 0.01$ (specifically, $R_3 = 0.0234$, $0.017$, and $0.0126$ for $\epsilon = 0.16$, $0.19$, and $0.22$, respectively). The fixed points at $R_1$ and $R_3$ are stable, as indicated by the negative slope of $G_{\epsilon,T}$ at those locations. These correspond to high and low activity states, respectively. In contrast, the intermediate fixed point at $R_2$ is unstable, since $G'_{\epsilon,T}(R_2) > 0$. (\textbf{b}) A saddle-node bifurcation occurs at $\epsilon \approx 0.24$, where the high-activity stable state and the intermediate unstable state coalesce and annihilate. The resulting bifurcation diagram, shown as a function of $\epsilon$, is constructed by numerically approximating the time-stepper in \eqref{eq: Time-stepper} using a Newton-Raphson solver. For $\epsilon = 0.24$, the system converges to a single stable steady state at $R \approx 0.0115$, corresponding to a low-activity regime.}
\label{fig:RHS1}
\end{figure}
\section*{Discussion}
In this study, we aimed to characterize the macroscopic dynamics of an agent-based network model using an  approach that requires significantly less information than traditional methods~\cite{Spil11}. 
For example, in our study, we use $50$ landmark agents for $R_t$, but $10000$ agents for $D_t$, which speeds up the analysis of the system. The investigation of macroscopic dynamics was accomplished by combining topological data analysis with the equation-free method, a framework we believe is applicable to a broad class of agent-based network models. Extending this methodology to other systems will be the focus of future work.

The most challenging component of our approach was addressing the inverse problem inherent in the implementation of the equation-free method, namely, reconstructing microscopic network configurations from macroscopic topological descriptors. The central topological feature we employed was the minimal value of the Rips filtration parameter (the minimum radius) at which the first persistent Betti number transitions from $0$ to $1$, signaling the emergence of a loop in the network's structure. This information served as a critical link between microscopic realizations and macroscopic behavior, enabling us to characterize and detect phase transitions in the network dynamics.

Within the equation-free framework, the inverse problem is addressed via a lifting operator, which typically relies on an optimization technique such as simulated annealing. While effective, simulated annealing is generally computationally intensive, particularly for large-scale agent-based networks, due to its reliance on stochastic sampling. To mitigate this inefficiency, we incorporated geometric information into the construction of the lifting operator. This modification significantly reduced the inherent randomness of the process, resulting in a more structured and computationally efficient lifting procedure.


Several other research works employ topological data analysis to examine dynamical systems \cite{topa15, Male16, xian2021capturing, Cole21, Guz22, Myer19}. For example, characterization of system dynamics and prediction of critical phase transitions have been conducted using CROCKER plots \cite{topa15}. To the best of our knowledge, no existing publication features CROCKER plots with an in-depth analysis of system dynamics, such as bifurcation and stability analysis. Our paper presents a detailed construction of the bifurcation diagram, integrating equation-free methods with topological data analysis. This approach paves the way for new, innovative directions combining CROCKER plots with the equation-free method. Our integration of the equation-free framework with topological data analysis suggests that equation-free methods can also be effectively combined with CROCKER plots to support more rigorous and quantitative analyses of complex systems. {A potential direction for future research is the application of our method to a broader class of agent-based models, including those of social networks, pedestrian dynamics, epidemic processes, and cell motility in biology and cancer. Such aggregation models are distinguished by emergent dynamics and spatially heterogeneous structures, rendering them a natural context for the application of our approach.}






%

\section*{Appendix}
\subsection*{Elements of Topological Data Analysis}
\label{sec:topology}

In this Supporting Material, we provide a brief introduction to Topological Data Analysis (TDA) and its application to complex dynamical systems. As these systems evolve, they generate data whose structure—or shape—can be examined and related to emergent system properties. TDA offers a framework for studying this shape by extracting topological features such as the number of connected components, circles, and voids. Through the filtration process, we analyze how these features appear and persist across scales, allowing us to identify the topological properties that remain robust.
 
We begin with the definition of the abstract simplicial complex and then proceed with the filtration process, finally, we describe the Betti numbers.


\subsubsection*{Simplicial complexes and simplicial homology}
Let $V=\{v_0,v_1,\dots,v_n\}$ be a finite set of points. A simplicial complex $K$ with vertex set $V$ is a collection of subsets of $V$ satisfying the following conditions:
\begin{itemize}
    \item All the singletons $\{v_j\}$ belong to $K$.
    \item Every subset of a set in $K$ belongs to $K$.
\end{itemize}
An element of $K$ containing $k+1$ vertices is called a $k$-simplex.

This definition is often referred to as a \emph{abstract simplicial complex}. Such a complex admits a canonical realization as a geometric simplicial complex by identifying each $0$-simplex $\{v_j\}$ with the vertex $v_j$ itself; each $1$-simplex $\{v_i,v_j\}$ with a line segment joining $v_i$ and $v_j$; each $2$-simplex $\{v_i,v_j,v_k\}$ with a triangle whose vertices are $v_i$, $v_j$, and $v_k$ and whose edges are identified with the line segments representing the $1$-simplices $\{v_i,v_j\}$, $\{v_i,v_k\}$, and $\{v_j,v_k\}$; and so on {Figure~\ref{fig:simplices}). This construction endows the simplicial complex with a canonical topology; see, for example, Ref.~\cite{Munkres_1984}.

\begin{figure}[H]
\centering
\includegraphics[width=1\linewidth]{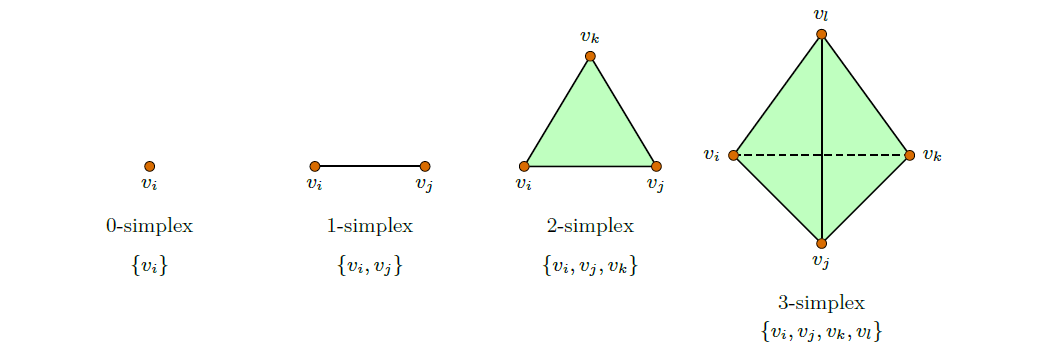}
\caption{Example of abstract simplicial complex for different cases.}
\label{fig:simplices}
\end{figure}
We define $C_k(K)$ to be the vector space over $\mathbb{Z}_2$ generated by the $k$-simplices of $K$. Equivalently, $C_k(K)$ consists of all formal linear combinations of $k$-simplices in $K$ with coefficients in the two-element field $\mathbb{Z}_2$. We further define linear maps $\partial_k \colon C_k(K) \to C_{k-1}(K)$, called boundary operators, by extending linearly over $\mathbb{Z}_2$ the following action on $k$-simplices:
\begin{equation}
    \partial_k\big(\{v_{j_0}, v_{j_1}, \dots, v_{j_k}\}\big)
    = \sum_{l=0}^k \{v_{j_0}, \dots, v_{j_{l-1}}, v_{j_{l+1}}, \dots, v_{j_k}\}.
\end{equation}
Elements of $C_k(K)$ are called $k$-chains. 

The chain spaces together with the boundary operators form a chain complex, that is, a sequence
\begin{equation}
    \{0\} \xrightarrow{\;\partial_{n+1}\;} C_n(K) \xrightarrow{\;\partial_n\;} \cdots \xrightarrow{\;\partial_{k+1}\;} C_k(K) \xrightarrow{\;\partial_k\;} C_{k-1}(K) \xrightarrow{\;\partial_{k-1}\;} \cdots \xrightarrow{\;\partial_2\;} C_1(K) \xrightarrow{\;\partial_1\;} \{0\}
\end{equation}
with the property that the composition of successive boundary operators vanishes, $\partial_k \partial_{k+1} = 0$. Elements of the kernel of $\partial_k$, that is, $k$-chains annihilated by $\partial_k$, are called $k$-cycles, while elements of the image of $\partial_{k+1}$ are called $k$-boundaries.
The $k$-cycles and $k$-boundaries form linear subspaces $Z_k(K)$ and $B_k(K)$, respectively, of the $k$-chain space $C_k(K)$, with $B_k(K) \subset Z_k(K)$ since every $k$-boundary is a $k$-cycle.
The $k$th homology group of $K$ (with coefficients in $\mathbb{Z}_2$) is defined as the quotient
\begin{equation}
    H_k(K) = Z_k(K) / B_k(K).
\end{equation}

Since the chain groups are vector spaces over $\mathbb{Z}_2$, the homology groups $H_k(K)$ are themselves $\mathbb{Z}_2$-vector spaces. Their dimensions,
\begin{equation}
    \beta_k(K) = \dim H_k(K),
\end{equation}
are called the Betti numbers of the simplicial complex $K$. Intuitively, $H_k(K)$ encodes the $(k+1)$-dimensional holes of $K$, and the $k$th Betti number $\beta_k(K)$ counts the number of such independent holes.

\subsubsection*{The Vietoris-Rips filtration}

\newcommand{\Rips}{\operatorname{VR}}
Suppose that $V$ is equipped with a distance function $d$. For each $r \geq 0$, the Vietoris--Rips simplicial complex at scale $r$ is defined by
\begin{equation}
    \{ v_{j_0}, v_{j_1}, \dots, v_{j_k} \} \in \Rips(V,r)
    \iff
    d(v_{j_a}, v_{j_b}) \leq r \text{ for all } 0 \leq a,b \leq k .
\end{equation}
Equivalently, a set of $k+1$ vertices spans a $k$-simplex in the Vietoris--Rips complex at scale $r$
if all pairwise distances are at most $r$. This construction yields a nested family of simplicial complexes, 
\begin{equation}
    r_1 \leq r_2 \implies \Rips(V,r_1) \subset \Rips(V,r_2),
\end{equation}
known as the Vietoris--Rips filtration.

As the scale parameter $r$ increases, additional simplices are incorporated into the complex, potentially altering its topological features, see Figure \ref{VRsimplex}.
These changes are quantified by the Betti numbers: the $k$th Betti number at scale $r$ is defined as the dimension of the $k$th homology group of $\Rips(V,r)$,
\begin{equation}
    \beta_k(\Rips(V,r)) = \dim H_k(\Rips(V,r)).
\end{equation}
Since $V$ is finite, the topology of the Vietoris--Rips complex can change only at finitely many values of the scale parameter $r$.

\begin{figure}[H]
\centering
\includegraphics[width=1.1\linewidth]{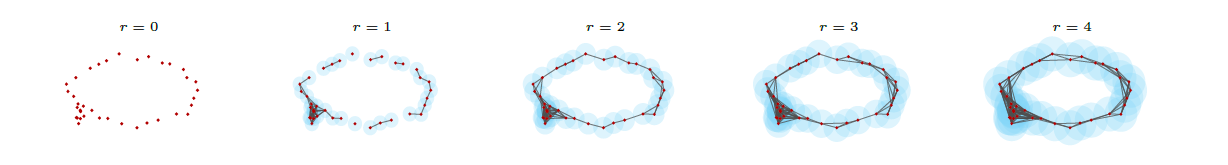}
\caption{Growth of the Vietoris--Rips complex as $r$ increases.}
\label{VRsimplex}
\end{figure}

\subsubsection*{Witness and Lazy Witness simplicial complexes}
For large datasets $X$, including every data point as a vertex in the Vietoris–Rips complex rapidly generates an excessively large set of simplices, making computations inefficient. The witness and lazy witness complexes address this challenge. These are sparse simplicial complexes designed to approximate the Vietoris–Rips complex of $X$ while significantly reducing computational cost for large point clouds.

In these constructions, a subset $L=\{l_0,l_1,\dots,l_n\}\subset X$, called the \emph{landmark points}, serves as the vertex set. The first landmark point $l_0$ is chosen randomly from $X$. Inductively, if $L_{i-1}=\{l_0,l_1,\dots,l_{i-2}\}$ denotes the set of the first $i-1$ landmark points, the $i$th landmark point is selected as the point in $X$ that maximizes the distance to $L_{i-1}$, where the distance from $x\in X$ to $L_{i-1}$ is defined as $d(x,L_{i-1})=\min\{d(x,l): l\in L_{i-1}\}$. 

The \emph{witness simplicial complex} $W(X,L,r)$ has vertex set $L$. A set of $k+1$ landmark points $\{l_{i_0},l_{i_1},\dots,l_{i_k}\}$ forms a $k$-simplex at scale $r$ if there exists a point $x\in X$, called a witness, such that
\begin{equation}
   \max\{d(l_{i_0}, x), d(l_{i_1},x),\dots, d(l_{i_k},x)\} \leq r + d(x,L).
\end{equation}
The witness complexes at different scales form a filtered simplicial complex, $W(X,L,r_1)\subseteq W(X,L,r_2)$ for $r_1\leq r_2$, providing an approximation of the Vietoris–Rips complex on $X$; see Figure \ref{fig:LWsimplex}.

The \emph{lazy witness complex} $LW(X,L,r)$ is a computationally simpler variant. It shares the same $1$-skeleton as the witness complex: the vertex set is $L$, and the $1$-simplices are defined using witnesses as
\begin{equation}
    \{l_{i_0},l_{i_1}\}\in LW(X,L,r) 
    \iff \max\{d(l_{i_0},x), d(l_{i_1},x)\} \leq r + d(x,L) \text{ for some } x\in X.
\end{equation}
A set of $k+1$ landmark points $\{l_{j_0},\dots,l_{j_k}\}$ then forms a $k$-simplex in $LW(X,L,r)$ if and only if all of its $2$-point subsets are $1$-simplices:
\begin{equation}
    \{ l_{j_0}, \dots, l_{j_k}\}\in LW(X,L,r) 
    \iff \{ l_{j_a}, l_{j_b}\}\in LW(X,L,r) \text{ for all } 0\leq a<b\leq k.
\end{equation}
The lazy witness complexes also form a filtered sequence of simplicial complexes that approximate the Vietoris–Rips complex on $X$; see Figure \ref{fig:LWsimplex}.
\begin{figure}[H]
\centering
\includegraphics[width=1.0\linewidth]{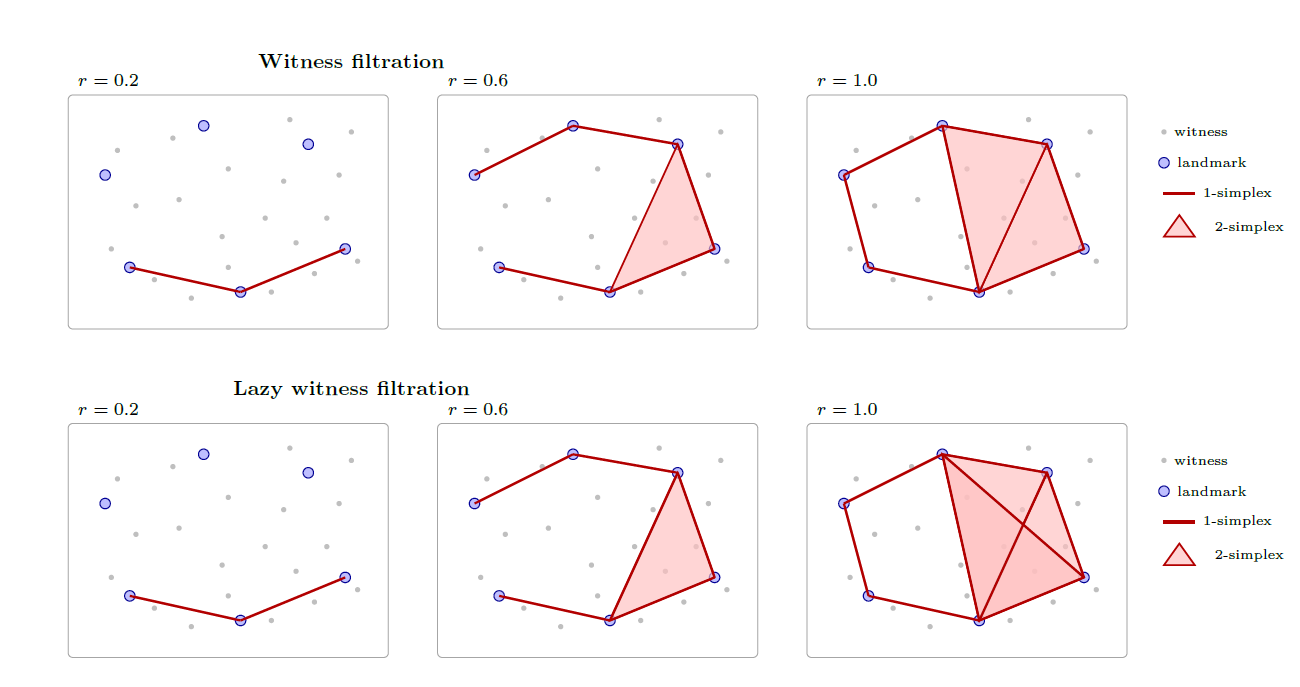}
\caption{Side-by-side visualization of witness vs. lazy witness filtrations on the same landmark set $L$ with a witness set $W$. Edit the edge/face lists to match your computed filtration.}
\label{fig:LWsimplex}
\end{figure}

\subsection*{Equation free method for complex systems}
\label{sec:equation_free_gener}
Consider a large-scale dynamical system defined on a network, where many interacting units collectively give rise to emergent behavior. Each unit—such as a neuron in a neural network, a cell in a biological tissue, or an individual in an epidemiological population—evolves in time according to well-specified microscopic rules. These rules describe the intrinsic dynamics of each node, which is transmitted through the connectivity. The interplay between the node-level dynamics and the network connectivity generates the global behavior of the system.

In general, we are interested in analyzing the network behavior on a different scale—the macroscopic one. Let $x \in \mathbb{R}^d$ ) denote the macroscopic variable, typically consisting of low-order statistics (e.g., density, mean value, variance) that characterize the emergent large-scale behavior of the system. Examples include the mean firing rate in neuronal networks or the magnetization in Ising-type models \cite{Kev09, Koz05, Spil11}.

Let $X \in \mathbb{R}^N$ be the microscopic state vector of the system, where each component describes the state of a single unit (e.g., one neuron), and where N $\gg d$ with $ d$  commonly equal to 1 or 2. The central assumption of the Equation-Free (EF) framework \cite{Gear03, Kev09, Sieb17, Marsch14} is that an effective macroscopic description of the system exists—typically in the form of ordinary differential equations, partial differential equations, or discrete-time maps—and that this macroscopic model closes in terms of the coarse variable $ x \in \mathbb{R}^d $. However, the explicit form of such an equation is not known, and in general it is extremely difficult to derive analytically, even under simplifying assumptions \cite{Kev09, Mont15}. For simplicity suppose that the macroscopic evolution law is described by the equation:
\begin{equation}
    x_{t+1}= F_T(x_t;\mu),
    \label{eq:dyn_sys_gen}
\end{equation}
where $x_t$ is the macroscopic variable $x$ on the macroscopic time $t$, and $\mu$ is the parameter of the system (for simplicity assume that $\mu \in \mathbb{R}$). The function $F_T:\mathbb{R}^d \rightarrow \mathbb{R}^d$ is not known (there is not an explicit formula), however will present a how we can numerically approximate at sampling time $T$. The macroscopic variable $x$ changes slowly in time relative to the microscopic vector state $X$. Starting with an initial value at time $t$ for the macroscopic variable $x$ i.e. $x(t)=x_t$, one has to transform this value to a microscopic network realization say $X_t \in \mathbb{R}^N$ i.e. to a state description for each node of the network,  consistent with $x_t$. This is done with a “lifting” operator, $L:\mathbb{R}^d \rightarrow \mathbb{R}^N$, such that
\begin{equation}
    X_t=L(x_t;\mu)
    \label{eq:lift_def}
\end{equation}
 This is a one-to-many operation (not unique), as there are in general many different network initialization that are consistent with $x$. The second step is to simulate the network rules (the evolution of neurons in the network)  for a short time, $T$ to obtain the new microscopic network state  $X_{t+1}$:
\begin{equation}
      X_{t+1} = \Phi_{T}(X_t;\mu)
\end{equation}
 The time $T$  is appropriate such that the macroscopic variable to change linear (similar to a tangent approximation). During the time $T$ the quickly-changing components of $X$ vary so that the probability density function of
$X$ becomes “slaved” to (or determined by) the current value of macrosopic $x$ \cite{laing06,Kev09,Gear03,Gear05}. Then, we define the restriction operator $R:\mathbb{R}^N \rightarrow \mathbb{R}^d$, to project back on the low dimensional space $\mathbb{R}^d$ where the macroscopic dynamics evolves:
\begin{equation}
    x_{t+1}=R(\Phi_T(X_t;\mu))
    \label{eq:Restr_def}
\end{equation}
This is a many-to-one operator, and results by calculating the first order statistics i.e. by averaging. The operators $L$ and $R$ satisfy the conditions $R\circ L = I $, the identity, so that the composition should have no effect, except an roundoff error. The function $F_T:\mathbb{R}^d \rightarrow \mathbb{R}^d $ is defined as the composition of the previous operators i.e.,
\begin{equation}
    F_T=R(\Phi_T(Lx_0;\mu)),
    \label{eq:F_t_def}
\end{equation}
Then, eq. \eqref{eq:dyn_sys_gen} is written
\begin{equation}
    x_{t+1}= F_T(x_t;\mu)=R(\Phi_T(Lx_t;\mu)),
    \label{eq:dyn_sys_gen_2}
\end{equation}
The aforementioned procedure constructs the function $F_T$ numerically, or the RHS of the macroscopic system evolution. In the case where the macroscopic description is one-dimensional i.e. $d=1$, since 
\begin{equation}
    F_T(x_t;\mu)=x_{t+1} 
    \label{eq:dyn_sys_gen_2_F}
\end{equation}
Then, $F_T$ can be directly constructed by discretizing the phase space and using the lift-evolve-restrict scheme, and the RHS is obtained.

 At this point, one can combine (utilise) computational algorithms like Newton’s or Newton-GMRES method to compute and trace branches of coarse-grained equilibria. Additionally, using iterative eigensolver methods for $F_T$, such as Arnoldi’s algorithm, can also be used to extract information about the stability of the coarse-grained system dynamics, overall to construct the bifurcation diagram \cite{SIET11,Gear03,Kev09}.

%
\printbibliography 
\end{document}